\documentclass[11pt, 3p]{elsarticle}

\usepackage{lineno}
\modulolinenumbers[5]

\usepackage[short]{optidef}
\usepackage{amsmath,amssymb,amsthm}
\usepackage{graphicx}
\usepackage{tabularx}
\usepackage{comment}
\usepackage[hidelinks]{hyperref}
\usepackage{adjustbox}
\usepackage{booktabs}
\usepackage{longtable}
\usepackage{multirow}
\usepackage{multicol}
\usepackage{enumitem}
\usepackage{array}
\usepackage{xcolor}
\usepackage{xcolor}

\usepackage{algorithm}
\usepackage{algpseudocode}
\usepackage{longtable}
\usepackage{multirow}
\usepackage{subcaption}

\journal{}

\newcommand{\MSSC}{\mathrm{MSSC}}
\newcommand{\cB}{\mathcal{B}}

\newcommand{\cP}{\mathcal{P}}
\newcommand{\cQ}{\mathcal{Q}}

\newtheorem{proposition}{Proposition}

\newtheorem{lemma}{Lemma}
\theoremstyle{remark}

\biboptions{authoryear}

\usepackage{tikz}
\usetikzlibrary{arrows.meta, shapes,calc,patterns,decorations.pathmorphing}
\usepackage{pgfplots}
\pgfplotsset{compat=1.16}
\usepackage{subfloat}

\tikzstyle{none}=[]
\tikzstyle{decision} = [diamond, aspect=2, draw, thick, text width=8em, text badly centered, inner sep=2pt]
\tikzstyle{decision2} = [diamond, aspect=3, draw, thick, text width=17em, text badly centered, inner sep=2pt]
\tikzstyle{io} = [trapezium, draw=black, thick, trapezium left angle=70, trapezium right angle=110, text width=17em, text centered]
\tikzstyle{start} = [rectangle, draw=black, thick, text width=17em, text centered, rounded corners, inner sep=5pt]
\tikzstyle{stop} = [rectangle, draw=black, thick, text width=5.5em, text centered, rounded corners, inner sep=5pt]
\tikzstyle{block} = [rectangle, draw=black, thick, text width=17em, text centered, inner sep=5pt]
\tikzstyle{block2} = [rectangle,
    text width=4em, text centered, inner sep=5pt]
\tikzstyle{line} = [draw, thick, -latex]

\usepackage{setspace}
\begin{document}

\begin{frontmatter}

\title{Constraint-Aware Anticlustering Decomposition
for Certifying Pairwise-Constrained Clustering}

\author{Anna Livia Croella}
\ead{annalivia.croella@unimercatorum.it}
\address{Department of Engineering and Science, Universitas Mercatorum, Rome, Italy}
\author{Veronica Piccialli}
\ead{veronica.piccialli@uniroma1.it}
\address{Department of Computer, Control and Management Engineering ``Antonio Ruberti'', Sapienza University of Rome, Rome, Italy}
\author{Antonio M. Sudoso}
\ead{antoniomaria.sudoso@unipd.it}
\address{Department of Mathematics ``Tullio Levi-Civita'', University of Padua, Padua, Italy}

\begin{abstract}
Certifying the quality of feasible solutions to large-scale minimum sum-of-squares clustering with pairwise must-link and cannot-link constraints remains challenging, as it requires strong lower bounds on the optimal objective value. We propose a decomposition framework that partitions the data points into smaller subsets, called anticlusters, and solves the induced clustering subproblems to obtain a certified optimality gap for the original problem. To preserve pairwise information, we introduce the Constrained Anticlustering Problem (CAP), which assigns each connected component of the constraint graph to a single anticluster while seeking geometrically representative subsets. We formulate CAP as a mixed-integer second-order cone program and develop a variable neighborhood search heuristic to construct high-quality decompositions. The resulting lower bound remains valid even when the induced subproblems are not solved to global optimality.
Extensive experiments on synthetic and real-world datasets with up to approximately \(4{,}000\) observations show that our approach can achieve mean optimality gaps below \(1\%\) at a modest computational cost.
\end{abstract}


\end{frontmatter}

\newpage

\section{Introduction}
Minimum sum-of-squares clustering (MSSC) is a fundamental problem in data analysis that has attracted extensive research in mathematical optimization and computer science. Given a set of observations and a prescribed number of clusters, MSSC seeks a partition that minimizes the sum of squared Euclidean distances from each observation to the centroid of its cluster. Classical $k$-means procedures approach this problem by alternating assignment and centroid updates \citep{MacQueen1967, HartiganWong1979, Lloyd1982}. Although computationally efficient, these procedures depend on initialization and provide no guarantee of global optimality.

Clustering is typically formulated as an unsupervised learning task, but applications often provide partial information about the desired partition. Pairwise must-link (ML) and cannot-link (CL) constraints offer a simple way to incorporate this knowledge: an ML constraint requires two observations to belong to the same cluster, whereas a CL constraint requires them to belong to different clusters. Because these relations do not require complete class labels, they are among the most widely studied forms of supervision in constrained clustering \citep{GonzalezAlmagroEtAl2025}.
The pairwise-constrained MSSC problem seeks a partition into a prescribed number of nonempty clusters that minimizes the within-cluster sum of squares while satisfying all ML and CL constraints. A mixed-integer nonlinear formulation represents this problem through binary assignment variables and continuous cluster centers.
Let $O=\{p_1,\ldots,p_N\}\subseteq\mathbb{R}^F$ be the set of observations and $K$ the prescribed number of clusters. For a positive integer $q$, write $[q]:=\{1,\ldots,q\}$. The pairwise information is specified by the sets of must-link pairs $\mathcal{ML}\subseteq[N]\times[N]$ and cannot-link pairs $\mathcal{CL}\subseteq[N]\times[N]$, with $\mathcal{B}: = \mathcal{ML} \cup\mathcal{CL}$.
For each $i\in[N]$ and $j\in[K]$, let $\delta_{ij}$ equal one if observation $p_i$ is assigned to cluster $j$ and zero otherwise, and let $\mu_j\in\mathbb{R}^F$ denote the center of cluster $j$. The problem can then be written as
\begin{subequations}
\label{def:MSSC_B}
\begin{align}
\operatorname{MSSC}_{\mathcal{B}}(O,K)
=
\min_{\substack{
\delta\in\{0,1\}^{N\times K}\\
\mu\in\mathbb{R}^{K\times F}}}
\quad
&\sum_{i=1}^N\sum_{j=1}^K
\delta_{ij}\|p_i-\mu_j\|^2_2
\label{def:MSSC_B_obj}
\\
\text{s.t.}\quad
&\sum_{j=1}^K\delta_{ij}=1,
&&\forall i\in[N],
\label{def:MSSC_B_assign}
\\
&\sum_{i=1}^N\delta_{ij}\ge 1,
&&\forall j\in[K],
\label{def:MSSC_B_nonempty}
\\
&\delta_{ij}=\delta_{i'j},
&&\forall(i,i')\in\mathcal{ML},\ \forall j\in[K],
\label{def:MSSC_B_ml}
\\
&\delta_{ij}+\delta_{i'j}\le 1,
&&\forall(i,i')\in\mathcal{CL},\ \forall j\in[K].
\label{def:MSSC_B_cl}
\end{align}
\end{subequations}
The objective~\eqref{def:MSSC_B_obj} minimizes the sum of squared Euclidean distances between observations and their assigned cluster centers. Constraints~\eqref{def:MSSC_B_assign} assign each observation to exactly one cluster, while constraints~\eqref{def:MSSC_B_nonempty} ensure that all clusters are nonempty. Constraints~\eqref{def:MSSC_B_ml} require the observations in each ML pair to share a cluster, whereas constraints~\eqref{def:MSSC_B_cl} require those in each CL pair to belong to different clusters. For any fixed feasible assignment, the optimal center $\mu_j$ is the centroid of the observations assigned to cluster $j$.

Pairwise information can make the resulting partition more consistent with application-specific knowledge, but it also further constrains an already difficult optimization problem. MSSC is NP-hard even without pairwise constraints \citep{AloiseEtAl2009, MahajanEtAl2012}. Heuristic methods can quickly produce feasible solutions to problem~\eqref{def:MSSC_B}, but their objective values provide only upper bounds on the global optimum. Certifying how close these solutions are to optimality also requires strong lower bounds. Exact methods provide such bounds, but their computational cost limits their applicability to large instances. Obtaining tight optimality certificates at this scale therefore remains challenging.

For large-scale unconstrained MSSC, \citet{CroellaPiccialliSudoso2027} proposed a decomposition approach that partitions the dataset into smaller, tractable subsets and sums the optimal values of the induced clustering subproblems to obtain a valid lower bound for the original problem. An auxiliary optimization problem guides the decomposition by constructing internally heterogeneous subsets, called \textit{anticlusters}, that reproduce the geometric structure of the full dataset.
Extending this approach to pairwise-constrained MSSC is not immediate since geometric representativeness alone is insufficient in this case: the decomposition must also preserve the pairwise constraints.
Indeed, if the points involved in a pairwise constraint are assigned to different anticlusters, that constraint is neglected. Although the resulting decomposition can still provide a valid lower bound, omitting constraints enlarges the subproblem feasible regions and can weaken the bound. 

\paragraph{Contributions.}
In this paper, we develop a decomposition framework for computing strong lower bounds and certified optimality gaps for MSSC with pairwise constraints.
We show that the optimality gap is given by the sum of two terms: the first arises
from decomposing the original problem into independent subproblems, whereas
the second results from using the available subproblem lower bounds in place
of their exact optimal values.
To preserve pairwise information, we introduce the \emph{Constrained Anticlustering Problem} (CAP). CAP represents the pairwise constraints through an undirected graph and assigns each connected component to a single anticluster, ensuring that every constraint is retained in exactly one induced subproblem. Its objective promotes geometric representativeness relative to a feasible reference clustering, while cardinality constraints control subproblem sizes. We formulate CAP as a mixed-integer second-order cone program (MISOCP) and develop a variable neighborhood search (VNS) heuristic to construct high-quality feasible partitions. Valid lower bounds for the induced clustering subproblems are then computed using a specialized exact solver. We evaluate the framework on synthetic and real-world constrained instances with up to $4{,}177$ observations. The experiments examine the trade-off between decomposition quality and subproblem tractability, compare VNS with direct MISOCP optimization, and assess the effect of preserving CL constraints. Extensive experiments show that the proposed framework achieves a mean optimality gap below \(1\%\), with mean total running times below \(10\) minutes, providing near-optimality guarantees for feasible solutions to large-scale constrained MSSC.

The remainder of the paper is organized as follows.
Section~\ref{sec:related_work} reviews the constrained clustering literature.
Section~\ref{sec:fundamentals} derives the decomposition bound for
pairwise-constrained MSSC.
Section~\ref{sec:problem_definition} presents CAP and the
lower bound computation.
Section~\ref{sec:computational} describes the experimental design and reports
the computational results.
Finally, Section~\ref{sec:conclusion} summarizes the main findings and
discusses directions for future research.
The Appendix provides detailed information on the benchmark
instances and computational experiments.

\section{Related work}
\label{sec:related_work}
Constrained clustering incorporates background knowledge into the clustering process through side constraints. Among the various forms of side information, instance-level must-link (ML) and cannot-link (CL) constraints are particularly common; a recent overview and taxonomy are provided by \citet{GonzalezAlmagroEtAl2025}. These constraints can be interpreted as hard constraints, which every feasible solution must satisfy, or as soft constraints, whose violations are incorporated into an objective or probabilistic model. 
Many constrained clustering methods adapt the alternating assignment and centroid-update scheme of unconstrained $k$-means. In the unconstrained setting, each observation can be assigned independently to its nearest center. Hard pairwise constraints couple these decisions, making feasibility itself a combinatorial issue alongside objective minimization. As discussed by \citet{DavidsonRavi2005}, a sequence of locally feasible assignments may fail to produce a complete feasible partition, even when one exists.
An influential early approach is COP-$k$-means \citep{WagstaffEtAl2001}, which assigns observations sequentially to nearby centers while respecting the prescribed hard constraints. Its greedy assignment procedure may stop without a feasible partition, even when a feasible partition exists. Related literature treats pairwise supervision through probabilistic models \citep{BasuEtAl2004}, metric learning \citep{BilenkoEtAl2004}, and penalty-based formulations, including constrained vector quantization \citep{DavidsonRavi2005, PellegBaras2007} and Lagrangian methods \citep{GanjiEtAl2016}. Metaheuristic approaches include dual iterative local search \citep{GonzalezAlmagroEtAl2020} and differential evolution \citep{GonzalezAlmagroEtAl2021}. 
These approaches establish important ways of using pairwise supervision, but changing the distance metric or allowing constraint violations leads to a different optimization model from the hard-constrained MSSC considered here.


Mathematical optimization can address hard constraints more explicitly. \citet{LeEtAl2019} formulate constrained $k$-means through binary optimization and develop an alternating scheme for searching the assignments. \citet{Baumann2020} proposes BLPKM, which solves a binary linear program for the assignment step with fixed centers and then updates the centers. This avoids the greedy construction used in COP-$k$-means. However, solving an assignment subproblem optimally does not establish global optimality of the joint problem.
The Pairwise-Confidence-Constraints-Clustering algorithm, PCCC, of \citet{BaumannHochbaum2025} handles hard and soft ML and CL constraints simultaneously and allows constraint-specific confidence weights. It combines integer programming assignments with ML preprocessing, an enhanced multistart strategy, and model-size reduction. Its reported experiments include instances with up to 60,000 observations and millions of CL constraints. Nevertheless, PCCC remains a heuristic: a feasible hard-constrained solution supplies an upper bound, whereas finding a valid lower bound to assess solution quality requires an additional bounding procedure.

Several optimization paradigms have been developed to solve constrained clustering globally. Column generation for ML- and CL-constrained MSSC is investigated by \citet{AloiseHansenRocha2012}. \citet{BabakiEtAl2014} develop a column generation framework that incorporates pairwise and other cluster constraints into the generation of feasible clusters. 
Constraint programming offers another way to combine clustering criteria with logical and combinatorial constraints. \citet{DaoEtAl2013} introduce a declarative framework for constrained clustering, subsequently developed in \citet{DaoEtAl2017}. For MSSC specifically, \citet{DaoEtAl2015} introduce a dedicated optimization constraint, a lower bound computation, and filtering procedures for the within-cluster sum-of-squares criterion. \citet{GunsEtAl2016} combine constraint programming with repetitive branch-and-bound, allowing user constraints to be integrated into the search. 
Semidefinite programming (SDP) provides particularly strong relaxations for MSSC \citep{PengWei2007, PiccialliSudosoWiegele2022}. For hard pairwise-constrained MSSC, \citet{PiccialliRussoSudoso2022} propose PC-SOS-SDP, a state-of-the-art solver that incorporates ML and CL constraints into an SDP-based branch-and-cut algorithm and strengthens its lower bounds through valid inequalities. Feasible clusterings are obtained using a rounding procedure based on the assignment step of \cite{Baumann2020}. Related SDP techniques also address prescribed cluster cardinalities \citep{PiccialliSudoso2025} and biclustering with ML and CL constraints \citep{sudoso2026exact}. These methods demonstrate the value of strong relaxations for certifying solution quality. Their computational cost, however, motivates decomposition approaches that derive bounds from smaller subproblems.

Subset-based bounds have a long history in exact MSSC algorithms. \citet{KoontzEtAl1975}, \citet{Diehr1985}, and \citet{Brusco2006} use bounds based on smaller clustering problems within branch-and-bound procedures. \citet{GunsEtAl2016} carry this subset-based bounding perspective into constrained MSSC. Thus, the general principle of obtaining bounds from smaller problems is the core of large-scale decomposition methods; the quality and computational usefulness of the resulting bound depend critically on how those problems are constructed. As shown in \cite{CroellaPiccialliSudoso2027}, anticlustering provides an effective way to guide this construction. In contrast to clustering, which groups similar observations, anticlustering creates internally heterogeneous subsets that are comparable to one another \citep{Spath1986, PapenbergKlau2021}. More recently, \citet{BaumannEtAl2026} propose the Assignment-Based Anticlustering algorithm for large Euclidean datasets.
For large-scale unconstrained MSSC, \citet{CroellaPiccialliSudoso2027} connect anticlustering directly to lower bound computation. 
The present work proposes an anticlustering-based bounding strategy for the constrained MSSC case.

\section{Decomposition bound for pairwise-constrained MSSC}
\label{sec:fundamentals}

Consider the pairwise-constrained MSSC problem defined in
\eqref{def:MSSC_B}.
Our approach seeks a partition
$\mathcal S=\{S_1,\ldots,S_T\}$ that decomposes the dataset into smaller
subsets. The constrained clustering problems induced on the subsets are used to construct a lower bound on $\textrm{MSSC}_{\mathcal{B}}(O, K)$.
The strength of this bound depends on the optimal values of the induced
subproblems: larger values yield a stronger bound. This motivates the
construction of internally heterogeneous subsets, called anticlusters,
that retain the structure of the complete dataset. A feasible reference clustering guides this construction.

For each $t\in[T]$, retain only the pairwise constraints whose endpoints
both belong to $S_t$. Specifically, define $\mathcal{ML}_t =\{(i,i')\in\mathcal{ML}:p_i,p_{i'}\in S_t\}$, $C_t =\{(i,i')\in\mathcal{CL}:p_i,p_{i'}\in S_t\}$
and let $\mathcal B_t:=\mathcal{ML}_t\cup\mathcal{CL}_t$.
Constraints whose endpoints belong to different anticlusters are
omitted and therefore relaxed in the decomposition.
The induced subproblem is $\operatorname{MSSC}_{\mathcal B_t}(S_t,K)$, with exactly $K$ nonempty clusters.
For any partition $\mathcal Q=\{A_1, \dots, A_q\}$ of $O$ into $q$ nonempty subsets, define its
within-cluster sum of squares by
\[
f(\mathcal Q)
:=
\sum_{A\in\mathcal Q}
\sum_{i:p_i\in A}\|p_i-\bar p_A\|^2_2,
\qquad
\bar p_A
:=
\frac{1}{|A|}\sum_{i:p_i\in A}p_i.
\]
This definition applies independently of the number of subsets in
$\mathcal Q$.

\begin{lemma}
\label{lem:restriction}
Let $\mathcal P=\{C_1,\ldots,C_K\}$ be a feasible partition for
$\operatorname{MSSC}_{\mathcal B}(O,K)$, and let
$\mathcal S=\{S_1,\ldots,S_T\}$ be a partition of $O$.
For each $t\in[T]$, define the restriction of $\mathcal{P}$ to $S_t$ as $\mathcal P|_{S_t}
:=
\{C_j\cap S_t:
j\in[K],\ C_j\cap S_t\neq\varnothing\}$.
If the induced subproblem
$\operatorname{MSSC}_{\mathcal B_t}(S_t,K)$ is feasible, then
\begin{equation}
\label{eq:fPSt}
\operatorname{MSSC}_{\mathcal B_t}(S_t,K)
\leq f(\mathcal P|_{S_t}).
\end{equation}
If the restriction contains exactly $K$ clusters, equality holds
if and only if it is optimal for the induced subproblem.
\end{lemma}

\begin{proof}
Assume that the induced subproblem $\MSSC_{\cB_t}(S_t,K)$ is feasible.
As noted after formulation~(1), for a fixed assignment, the centroid
minimizes the sum of squared distances within each cluster. Thus
$\MSSC_{\cB_t}(S_t,K)$ is the minimum of $f(\cQ)$ over all partitions
$\cQ$ of $S_t$ into exactly $K$ nonempty clusters.
We first observe that $\cP|_{S_t}$ satisfies every constraint in $\cB_t$. Indeed, if
$(i,i') \in \mathcal{ML}_t$, then $p_i, p_{i'} \in C_j \cap S_t$ for some $j$, whereas if
$(i,i') \in \mathcal{CL}_t$, then $p_i$ and $p_{i'}$ lie in $C_j \cap S_t$ and
$C_{j'} \cap S_t$ with $j \neq j'$.
The restriction may, however, contain fewer than $K$ nonempty clusters.
Let $\Gamma_1,\dots,\Gamma_m$ be the connected components of the graph with vertex set $S_t$ and edge set
$\{\{p_i,p_{i'}\} : (i,i') \in \mathcal{ML}_t\}$. All observations in a component
must belong to the same cluster in any partition satisfying $\mathcal{ML}_t$.
Thus every cluster of $\cP|_{S_t}$ is
a union of components. Since a feasible partition has $K$ nonempty
clusters, we have $m \ge K$.
Let $q$ be the number of clusters in $\cP|_{S_t}$. Whenever the current
partition has fewer than $K$ clusters, some cluster $A$ contains at
least two components. Choose one such component $\Gamma$ and split
$A$ into the nonempty clusters $A_1 = \Gamma$ and
$A_2 = A \setminus \Gamma$. Each split keeps every must-link component
intact and preserves the separation of every cannot-link pair.
After $K-q$ splits, we obtain a feasible partition $\cQ$ of $S_t$ into $K$ nonempty clusters that refines
$\cP|_{S_t}$.
Each split does not increase the objective. Indeed, using the centroid
$\bar p_{A_\ell}$ for each new cluster $A_\ell$ gives
\begin{equation*}
  \sum_{i:\,p_i \in A} \|p_i - \bar p_A\|^2_2
  = \sum_{\ell=1}^2 \sum_{i:\,p_i \in A_\ell}
     \|p_i - \bar p_A\|^2_2 \ge \sum_{\ell=1}^2 \sum_{i:\,p_i \in A_\ell}
     \|p_i - \bar p_{A_\ell}\|^2_2.
\end{equation*}
All other clusters retain their contributions to the objective. Applying
this inequality at each split yields $\MSSC_{\cB_t}(S_t,K) \le f(\cQ) \le f(\cP|_{S_t})$, which proves \eqref{eq:fPSt}. If $q=K$, then $\cP|_{S_t}$ is itself feasible, and equality in \eqref{eq:fPSt} holds if and only if
$\cP|_{S_t}$ is optimal for the induced subproblem.
\end{proof}

We next recall the within-cluster sum-of-squares decomposition used
in \citet{CroellaPiccialliSudoso2027}. The following geometric identity does
not depend on the presence of pairwise constraints.

\begin{lemma}
\label{lemma:anova}
Let $\mathcal{P}=\{C_1,\dots,C_K\}$ be a partition of $O$, and let
$\mu_j:=\frac{1}{|C_j|}\sum_{i:p_i\in C_j}p_i$ be the centroid of $C_j$. For a partition $\mathcal{S}=\{S_1,\dots,S_T\}$ of $O$, define
$A_{jt}:=C_j\cap S_t$ and $n_{jt}:=|A_{jt}|$. Whenever $n_{jt}>0$, let
$\mu_{jt}:=\frac 1 {n_{jt}}\sum_{i:p_i\in A_{jt}}p_i$ be the centroid of $A_{jt}$.
Then, for every $j\in[K]$,
\begin{equation}
\label{eq:anova}
\sum_{i:p_i\in C_j}\|p_i-\mu_j\|^2_2
=
\sum_{t \in [T] : n_{jt}>0} \sum_{i:p_i\in A_{jt}}\|p_i-\mu_{jt}\|^2_2
+
\sum_{t \in [T] : n_{jt}>0} n_{jt}\|\mu_{jt}-\mu_j\|^2_2.
\end{equation} 
\end{lemma}
Summing \eqref{eq:anova} over the clusters yields
\begin{equation}
\label{eq:reference_decomposition}
f(\mathcal P)
=
\sum_{t=1}^T f(\mathcal P|_{S_t})
+
\sum_{j=1}^K
\sum_{t \in [T] : n_{jt}>0}
n_{jt}\|\mu_{jt}-\mu_j\|^2_2.
\end{equation}

Lemma \ref{lemma:anova} shows that the contribution of a cluster $C_j$ can be decomposed into two parts.
The first is the sum of the within-cluster objective values associated
with the intersections $A_{jt}=C_j\cap S_t$.
The second is a term measuring the deviation of the induced subproblem centroids from $\mu_j$.
For a fixed reference partition, smaller centroid deviations therefore
preserve a larger portion of its objective value in the restricted
partitions. However, the induced subproblems may improve on these
restrictions, decreasing the bound. Centroid preservation thus motivates the construction
of anticlusters but does not by itself guarantee a strong lower bound.

We now combine these results to obtain a lower bound on the
optimal value of problem \eqref{def:MSSC_B}.
To account for induced subproblems that may be infeasible, we adopt the convention that their contribution to the bound is zero.
Specifically, we define
\begin{equation*}
W_t
:=
\begin{cases}
\operatorname{MSSC}_{\mathcal B_t}(S_t,K),
& \text{if the induced subproblem is feasible},\\
0,
& \text{otherwise}.
\end{cases}
\end{equation*}
For feasible subproblems, Lemma~\ref{lem:restriction} ensures that
$W_t\leq f(\mathcal P|_{S_t})$ for every feasible
partition $\mathcal P$. For infeasible subproblems, the same inequality
follows from the nonnegativity of the within-cluster sum of squares objective.

\begin{proposition}
\label{prop:lower_cl}
For every partition $\mathcal S=\{S_1,\ldots,S_T\}$ of $O$,
\begin{equation}
\label{eq:decomposition_lb}
\operatorname{MSSC}_{\mathcal B}(O,K)
\geq
\sum_{t=1}^T W_t
=:LB^*(\mathcal S).
\end{equation}
If every induced subproblem is feasible, then
$
LB^*(\mathcal S)
=
\sum_{t=1}^T
\operatorname{MSSC}_{\mathcal B_t}(S_t,K).
$
\end{proposition}

\begin{proof}
Let $\mathcal P^*=\{C_1^*,\ldots,C_K^*\}$ be an optimal
partition for $\operatorname{MSSC}_{\mathcal B}(O,K)$.
For each $j\in[K]$ and $t\in[T]$, let
$n_{jt}^*:=|C_j^*\cap S_t|$, and denote the corresponding global
and local centroids by $\mu_j^*$ and $\mu_{jt}^*$, respectively,
with local centroids defined only for nonempty intersections.
Applying \eqref{eq:reference_decomposition} gives
\begin{equation}
\label{eq:decomposition}
\operatorname{MSSC}_{\mathcal B}(O,K)
=
\sum_{t=1}^T f(\mathcal P^*|_{S_t})
+
\sum_{j=1}^K
\sum_{t\in[T] : n_{jt}^*>0}
n_{jt}^*\|\mu_{jt}^*-\mu_j^*\|^2_2.
\end{equation}
The second term is nonnegative. Combining this observation with
$W_t\leq f(\mathcal P^*|_{S_t})$ for every $t\in[T]$ yields
\[
\operatorname{MSSC}_{\mathcal B}(O,K)
\geq
\sum_{t=1}^T f(\mathcal P^*|_{S_t})
\geq
\sum_{t=1}^T W_t,
\]
which proves \eqref{eq:decomposition_lb}.
\end{proof}

In practice, the induced subproblems need not be solved to optimality.
For each feasible subproblem, its optimal value can be replaced by
any valid lower bound $LB_t\leq W_t$; infeasible subproblems
contribute zero by convention. The sum of these contributions remains
a valid lower bound on the complete instance. 

Proposition~\ref{prop:lower_cl} establishes the validity of the
decomposition bound but does not specify how to select a strong
anticlustering partition. Ideally, one would maximize
$LB^*(\mathcal S)$ directly. However, evaluating this quantity for
each candidate partition requires solving all induced constrained
clustering problems. The next section introduces a construction guided by a feasible reference clustering and the pairwise-constraint structure.

\section{The Constrained Anticlustering Problem}
\label{sec:problem_definition}

We use a feasible clustering $\mathcal P$ as a reference and seek anticlusters satisfying three requirements: (i) they should reproduce the geometric structure of $\mathcal P$, (ii) preserve every pairwise constraint,
and (iii) contain a balanced representation of the reference clusters.
The first requirement is captured by minimizing the deviations from the reference centroids so that each anticluster represents the distribution of the complete dataset.
The second prevents the induced subproblems from relaxing pairwise constraints.
The third controls their size and composition, producing subproblems with comparable numbers of observations from each reference cluster. Together, these requirements define the \emph{Constrained Anticlustering Problem} (CAP).

The reference $\mathcal{P}$ has two distinct roles. Its objective value
$UB=f(\mathcal P)$ provides an upper bound for the complete constrained MSSC problem, while its assignments and centroids guide the construction of $\mathcal S$. The quality of $\mathcal P$ may affect the anticlustering partition selected by CAP, but not the validity of the lower bound computed from that partition.
For a candidate anticlustering partition $\mathcal{S}$, we define the optimality gap as
\begin{equation}\label{eq:gapstar}
    \gamma^*(\mathcal{S}):=
    \frac{UB-LB^*(\mathcal{S})}{UB}.
\end{equation}

\begin{figure}[ht]
    \centering
    \includegraphics[width=\textwidth]
    {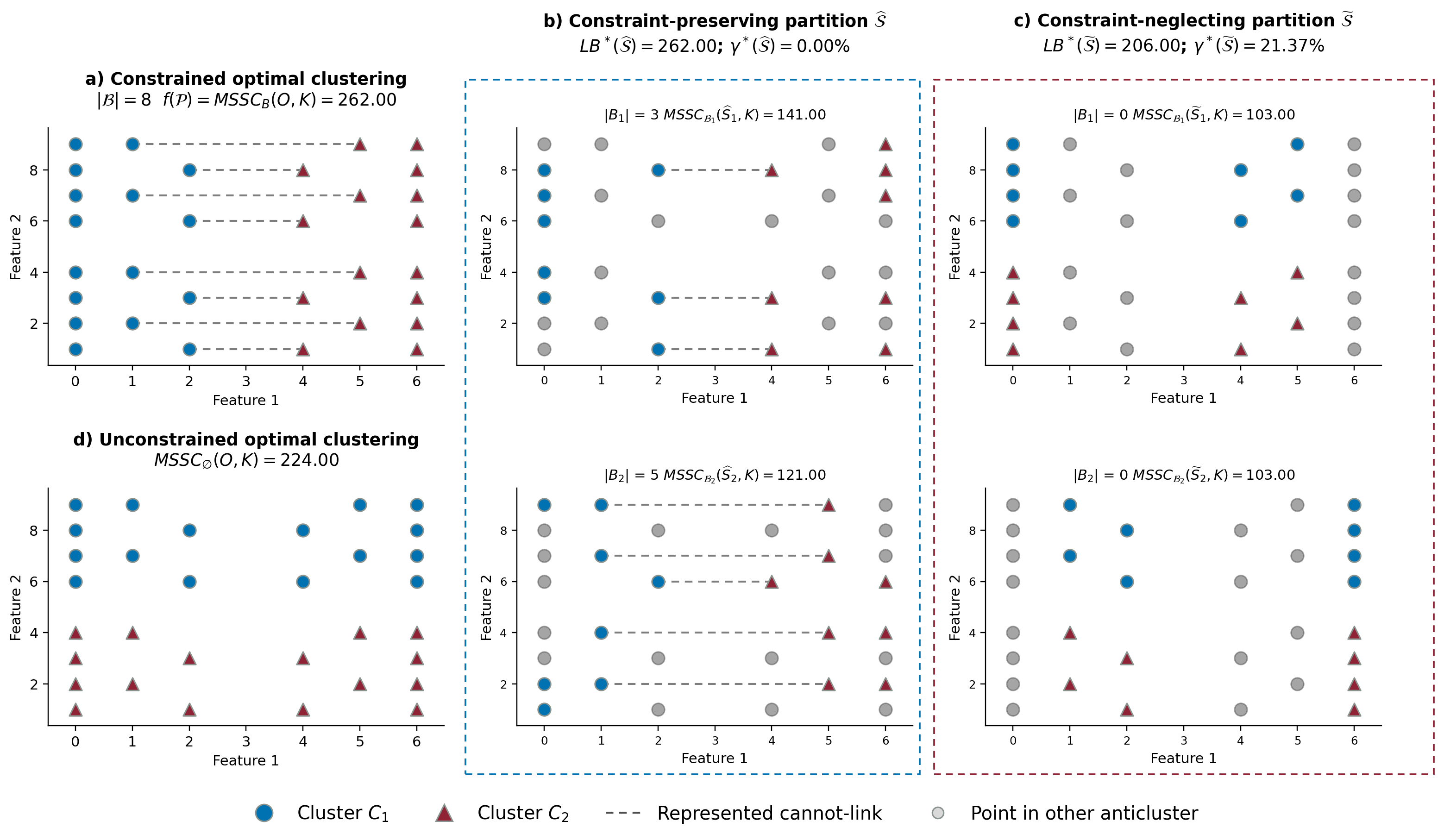}
    \caption{Effect of the anticlustering partition on the decomposition bound for an instance with $K=2$ clusters and $T=2$ anticlusters.}
    \label{fig:CAP_example}
\end{figure}
Figure~\ref{fig:CAP_example} illustrates the roles of geometric representativeness and constraint preservation on an instance with $K=2$ clusters and $T=2$ anticlusters.
Panel~(a) shows an optimal clustering of the complete constrained instance, with objective value
$\operatorname{MSSC}_{\mathcal B}(O,2)=262$. The two marker types identify the clusters, while the dashed segments represent the cannot-link constraints. Panels~(b) and~(c) show two alternative partitions $\hat{\mathcal S}=\{\hat S_1,\hat S_2\}$ and
$\tilde{\mathcal S}=\{\tilde S_1,\tilde S_2\}$. In each subpanel, colored points belong to the anticluster under consideration and retain the marker of their reference cluster, whereas gray points belong to the other anticluster. Finally, panel~(d) shows the optimal clustering of the same instance after removing the pairwise constraints, whose objective value is $\operatorname{MSSC}_{\varnothing}(O,2)=224$.
The partition $\hat{S}$ in panel~(b) preserves all cannot-link constraints. The two subproblem values are
$141$ and $121$, and hence $LB^*(\mathcal{\hat{S}})=262$. In this example,
the centroid of each intersection $C_j\cap \hat{S}_t$ coincides with the centroid of
$C_j$ given by the clustering $\mathcal{P}$. The nonnegative centroid-deviation term in
Lemma~\ref{lemma:anova} is therefore zero, and the lower bound equals the optimal
value of the complete constrained instance, i.e. $\gamma^*(\mathcal{\hat{S}}) = 0$.
The partition $\tilde{S}$ in panel~(c), instead, neglects all ten cannot-link pairs. In this case, the two subproblem values are $103$ and $103$, giving
$LB^*(\mathcal{\tilde{S}})=206$. The corresponding optimality gap is $\gamma^*(\mathcal{\tilde{S}}) = 21.37\%$.
Note that, once a cannot-link pair is neglected,
no induced subproblem enforces the separation of its endpoints. The
subproblem solutions may consequently move toward the unconstrained structure shown in panel~(d); if an induced subproblem retains no cannot-link constraint,
it is itself an unconstrained MSSC problem.
Both partitions yield valid lower bounds, but the example shows that constraint preservation can be decisive
for their strength.

\subsection{Mixed-integer second-order cone formulation}
\label{sec:ac_construction}

To construct anticlusters that preserve all pairwise constraints, we first represent the constraint structure through the undirected graph
$G_{\mathcal B}=(V,E)$, where $V=[N]$ and
$E=\{\{i,i'\}:(i,i')\in\mathcal B\}$. Each edge therefore represents either an ML or a CL constraint. Let $U_1,\ldots,U_R$ denote the connected components of $G_{\mathcal B}$, with isolated vertices treated as singleton components.
The collection $\{U_1,\ldots,U_R\}$ forms a partition of $[N]$. Since the endpoints of every pairwise constraint belong to the same connected component, assigning these components as indivisible units ensures that no constraint is split across different anticlusters.
Importantly, this applies only when constructing the anticlustering partition. In the clustering problem, the original semantics of the constraints remain unchanged: the endpoints of an ML pair must belong to the same cluster, whereas those of a CL pair must belong to different clusters. Consequently, a connected component of $G_{\mathcal B}$ may contain observations assigned to different clusters in the reference $K$-partition.

Given a feasible $K$-partition $\mathcal{P}$, summing identity~\eqref{eq:anova} over $K$ clusters gives
\begin{equation*}
\label{eq:dsum}
f(\mathcal P)
=
\sum_{t=1}^T f(\mathcal P|_{S_t})
+
D(\mathcal P,\mathcal S),
\end{equation*}
where the aggregate centroid deviation is defined as
\begin{equation*}
\label{eq:def_D}
D(\mathcal P,\mathcal S)
:=
\sum_{j=1}^K
\sum_{t \in [T] : n_{jt}>0}
n_{jt}\|\mu_{jt}-\mu_j\|^2_2.
\end{equation*}
Since the reference partition is fixed, $f(\mathcal P)$ is constant.
Minimizing $D(\mathcal P,\mathcal S)$ is therefore equivalent to maximizing
the sum of the objective values induced by the fixed assignments on the
anticlusters, $\sum_{t=1}^T f(\mathcal P|_{S_t})$. By
Lemma~\ref{lem:restriction}, each of these values is an upper bound on the
optimal value of the corresponding induced subproblem. Hence,
$D(\mathcal P,\mathcal S)$ provides a tractable geometry-based surrogate for
selecting an anticlustering partition with a strong decomposition bound.

Accordingly, CAP assigns the connected components
$U_1,\ldots,U_R$ to the $T$ anticlusters while minimizing
$D(\mathcal P,\mathcal S)$ and balancing the representation of the reference clusters. The formulation must therefore express the centroid deviation and the balancing requirements in terms of component-to-anticluster assignments.
We next introduce the component-level parameters and decision variables
needed to represent these requirements. The resulting optimization problem is a mixed-integer second-order cone program.

For each connected component $U_r$ and reference cluster $C_j$, let
$a_{rj}:=|\{i\in U_r:p_i\in C_j\}|$. These counts are needed because a connected component may intersect more than one reference cluster.
Let $x_{rt}$ be a binary variable that equals one if the connected component $U_r$ is assigned to anticluster $S_t$, and zero otherwise.
We also introduce continuous variables $n_{jt}\ge0$ and $z_{jt}\ge0$, where $n_{jt}$ counts the points of $C_j$ assigned to $S_t$ and $z_{jt}$ represents their contribution to $D(\mathcal{P},\mathcal{S})$. Define $v_{jt}(x)
:=
\sum_{r=1}^{R}
\sum_{\substack{i\in U_r:\, p_i\in C_j}}
(p_i-\mu_j)x_{rt}$ as the sum of the differences between the centroid of cluster $j$ and the points assigned to the anticluster $t$ as the sum of the deviations from $\mu_j$ of the points of $C_j$ assigned to $S_t$.
The CAP formulation reads as follows
\begin{subequations}
\label{model:component_ac_miqcp}
\begin{align}
\min_{\substack{
x\in\{0,1\}^{R\times T}\\
n, z \in \mathbb{R}_+^{K \times T}}}
\quad
&\sum_{j=1}^K\sum_{t=1}^T z_{jt}
\label{model:component_ac_obj}
\\
\text{s.t.}\quad
&\sum_{t=1}^T x_{rt}=1,
&&\forall r\in[R],
\label{model:component_ac_assignment}
\\
&n_{jt}=\sum_{r=1}^R a_{rj}x_{rt},
&&\forall j\in[K],\ \forall t\in[T],
\label{model:component_ac_count}
\\
&\ell_j(T)\le n_{jt}\le u_j(T),
&&\forall j\in[K],\ \forall t\in[T],
\label{model:component_ac_cluster_balance}
\\
&\lVert v_{jt}(x)\rVert_2^2
\leq z_{jt}n_{jt},
&&\forall j\in[K],\ \forall t\in[T],
\label{model:component_ac_cone}
\end{align}
\end{subequations}
Constraint~\eqref{model:component_ac_assignment} assigns each connected
component, and hence all its points, to exactly one anticluster.
Equalities~\eqref{model:component_ac_count} define the number \(n_{jt}\)
of points from cluster \(C_j\) assigned to anticluster \(S_t\).
Although \(n_{jt}\) is declared continuous, it is integer-valued whenever
\(x\) is binary, since the coefficients \(a_{rj}\) are integer.
Constraints~\eqref{model:component_ac_cluster_balance} impose the required
balance for every cluster and anticluster.
Since \(n_{jt}\geq 0\) and \(z_{jt}\geq 0\), the inequality \eqref{model:component_ac_cone} is a rotated second-order cone constraint. Equivalently, it can be
written as the standard second-order cone constraint
\begin{equation}
\left\|
\begin{bmatrix}
2v_{jt}(x)\\
z_{jt}-n_{jt}
\end{bmatrix}
\right\|_2
\leq z_{jt}+n_{jt}.
\end{equation}
For an integer feasible solution with \(n_{jt}>0\), we have
$    v_{jt}(x) =
\sum_{i:\,p_i\in C_j\cap S_t}(p_i-\mu_j) =
n_{jt}(\mu_{jt}-\mu_j)$.
Therefore, constraint~\eqref{model:component_ac_cone} implies
$z_{jt}
\geq
\lVert v_{jt}(x)\rVert_2^2 / n_{jt}
=
n_{jt}\lVert\mu_{jt}-\mu_j\rVert_2^2$.
Since \(z_{jt}\) is minimized, this inequality is tight at optimality.
If \(n_{jt}=0\), then \(v_{jt}(x)=0\), and minimization yields
\(z_{jt}=0\). Thus, the objective of problem \eqref{model:component_ac_miqcp} is exactly
\(D(\mathcal{P},\mathcal{S})\).
The continuous relaxation is a convex second-order cone program.
Nevertheless, the binary assignment variables make the problem
combinatorial, and solving it to optimality can be computationally
expensive for large instances. We therefore use a heuristic procedure to obtain high-quality feasible solutions.

Note that, since connected components cannot be split, the balance bounds constrain the existence of a feasible CAP solution: some choices of $T$, $\ell_j(T)$, and $u_j(T)$ may render formulation~\eqref{model:component_ac_miqcp} infeasible.
Section~\ref{subsec:settings} specifies the balance bounds used in the
experiments and explains how we retain only values of $T$ for which
a feasible assignment exists.

\subsection{Variable neighborhood search heuristic}
\label{sec:vns}
In this section, we describe the variable neighborhood search (VNS) heuristic
\citep{mladenovic1997vns} for problem \eqref{model:component_ac_miqcp}.
The contribution of $U_r$ to reference cluster $C_j$ is $b_{rj}:=\sum_{i\in U_r:\, p_i\in C_j}(p_i-\mu_j)$.
Let $\sigma=(\sigma_1,\ldots,\sigma_R)$ denote an assignment of components to the anticlusters.
The VNS maintains the counts
$n_{jt}(\sigma):=\sum_{r:\sigma_r=t}a_{rj}$ and centered sums
$v_{jt}(\sigma):=\sum_{r:\sigma_r=t}b_{rj}$ for each cluster--anticluster
pair $(j,t)$.
Its objective value is evaluated as
\begin{equation}
\label{eq:vns_objective}
H(\sigma):=
\sum_{j=1}^K\sum_{t=1}^T
\begin{cases}
\dfrac{\|v_{jt}(\sigma)\|^2_2}{n_{jt}(\sigma)},& n_{jt}(\sigma)>0,\\
0,&n_{jt}(\sigma)=0.
\end{cases}
\end{equation}
Thus, $H(\sigma)$ coincides with $D(\mathcal{P},\mathcal{S})$ for the partition induced by $\sigma$.
Algorithm~\ref{alg:vns} summarizes the procedure. The heuristic uses a multistart strategy with
$n_{\mathrm{start}}$ starts, each initialized from a randomly generated feasible assignment of connected components to anticlusters. At each start, VNS improves the initial assignment through local search and shaking. The best solution across all starts is retained, and its assignment defines a feasible CAP solution.

\begin{algorithm}[!ht]
\caption{VNS heuristic for CAP}
\label{alg:vns}
\footnotesize
\begin{algorithmic}[1]
\setlength{\itemsep}{-1.2pt}
\setlength{\parsep}{0pt}
\setlength{\topsep}{0pt}
\setlength{\partopsep}{0pt}
\Require Connected components $U_1,\ldots,U_R$, aggregates $(a_{rj},b_{rj})$,
balance bounds $(\ell_j,u_j)$, number of starts $n_{\mathrm{start}}$, maximum number of attempts
$\alpha_\textrm{init}$,
maximum number of VNS iterations $q_{\max}$, and maximum neighborhood size
$\kappa_{\max}$
\Ensure A feasible anticlustering assignment $\bar{x}$
\State $\sigma^\star\gets\varnothing$, $H^\star\gets+\infty$
\For{$s=1,\ldots,n_{\mathrm{start}}$}
    \State $\sigma\gets\textsc{RandomAssignment}(U,a,\ell,u, \alpha_\textrm{init})$
    \State $\sigma\gets\textsc{LocalSearch}(\sigma)$ 
    \State $\kappa\gets1$, $q\gets0$
    \While{$q<q_{\max}$}
        \State $\sigma'\gets\textsc{Shake}(\sigma,\kappa)$ 
        \State $\sigma'\gets\textsc{LocalSearch}(\sigma')$
        \If{$H(\sigma')<H(\sigma)$}
            \State $\sigma\gets\sigma'$, $\kappa\gets1$
        \Else
            \State $\kappa\gets1+(\kappa\bmod\kappa_{\max})$
        \EndIf
        \State $q\gets q+1$
    \EndWhile
    \If{$H(\sigma)<H^\star$}
        \State $\sigma^\star\gets\sigma$, $H^\star\gets H(\sigma)$
    \EndIf
\EndFor
\State $\bar{x}_{rt}\gets1$ if $\sigma_r^\star=t$, and $\bar{x}_{rt}\gets0$ otherwise
\State \Return $\bar{x}$
\end{algorithmic}
\end{algorithm}

The routine \textsc{RandomAssignment} uses randomized
backtracking. Connected components are processed in nonincreasing
order of cardinality, with ties broken by the nonincreasing number
of reference clusters they intersect and then by component index.
At each construction step, candidate anticlusters satisfying the
cluster-wise upper bounds are considered in random order.
If no candidate remains, the procedure backtracks to a previous
component. An assignment is accepted only
if it also satisfies all lower bounds; otherwise, a new attempt
is made. If no feasible assignment is found within the prescribed
number of attempts $\alpha_\textrm{init}$, the corresponding start is skipped.

The routine \textsc{LocalSearch} is a first-improvement descent with two neighborhoods:
relocation of one connected component and exchange of two connected
components assigned to different anticlusters. Component order is
randomized at each descent iteration. Relocations are examined first;
exchanges are considered only if no improving relocation is found.
After any accepted move, the search restarts with the relocation
neighborhood. All moves preserve the cluster-wise lower and upper
bounds, and a move is accepted only if it decreases $H$. Incremental updates to the counts and centered
sums allow move evaluation without reconstructing the partition.

The routine \textsc{Shake} applies $\kappa$ random feasible moves. At each shaking step, a relocation
or an exchange is selected with equal probability.
Local search is then applied to the perturbed assignment. A solution improving the incumbent is accepted, and $\kappa$ is reset to one; otherwise, $\kappa$ is increased by one and reset to one after $\kappa_{\max}$. Each start terminates after at most $q_{\max}$ shaking--local-search
iterations.

The VNS and the quality of $\mathcal{P}$ can therefore affect the strength of the resulting bound, but not its validity, which follows from Proposition~\ref{prop:lower_cl}.

\subsection{Lower bound computation}
\label{sec:pc-sos-lb}
After choosing $\mathcal{S}$, we handle the constrained MSSC problem induced by each $S_t$ with an exact solver, 
which (i) proves the subproblem optimum, (ii) returns a valid bound
$LB_t$ or (iii) proves that the subproblem is infeasible.
Summing these bounds over $T$ and
applying Proposition~\ref{prop:lower_cl} gives the complete chain
\begin{equation}
\label{eq:bound_chain}
LB:=\sum_{t=1}^T LB_t
\le
LB^*(\mathcal{S})
\le
\operatorname{MSSC}_{\mathcal{B}}(O,K)
\le UB.
\end{equation}

If every subproblem is solved to optimality, then $LB=LB^*(\mathcal{S})$ and
$\gamma^*(\mathcal{S})$ are the optimality gaps defined
in \eqref{eq:gapstar}. If some subproblems terminate before optimality or are infeasible, their available lower
bounds still give a valid, possibly larger gap defined as:
\begin{equation}
\label{eq:gamma_lb}
\widehat\gamma (\mathcal{S}):=\frac{UB-LB}{UB}=
\frac{UB-LB^*(\mathcal{S})}{UB}
+
\frac{LB^*(\mathcal{S})-LB}{UB} = \gamma^{\mathrm{part}}(\mathcal{S}) + \gamma^{\mathrm{solve}}(\mathcal{S}).
\end{equation}
We call $\gamma^{\mathrm{part}}(\mathcal{S})$ the \emph{partition gap} and
$\gamma^{\mathrm{solve}}(\mathcal{S})$ the \emph{subproblem-solution gap}.
The partition gap is the residual gap with respect to the feasible reference
value $UB$ that would remain if all induced subproblems were solved to
optimality. It therefore measures the loss associated with the selected
anticlustering partition, relative to that reference solution. The
subproblem-solution gap measures the loss due to terminating the subproblem
solves before their optimal values are known; it is zero whenever every
subproblem is solved to global optimality.

Assume that all induced subproblems
$\textrm{MSSC}_{\mathcal{B}_t}(S_t,K)$ are feasible.
Let $\widetilde{UB}_t$ be an upper bound for subproblem $t$ and define
$
UB_t:=\min\bigl\{\widetilde{UB}_t,f(\mathcal{P}|_{S_t})\bigr\}
$ and $
LB^+:=\sum_{t=1}^T UB_t.
$
Since $UB_t$ does not exceed the value of the restricted reference solution, the available quantities satisfy
$LB\le LB^*(\mathcal{S})\le LB^+\le UB.$
Although they do not identify the two terms in
\eqref{eq:gamma_lb} exactly, we have
\begin{equation}
\label{eq:observable_gap_bounds}
\underline{\gamma}^{\mathrm{part}}(\mathcal{S}) = \frac{UB-LB^+}{UB}
\le
\gamma^{\mathrm{part}}(\mathcal{S}),
\qquad
\gamma^{\mathrm{solve}}(\mathcal{S})
\le
\frac{LB^+-LB}{UB} =
\overline{\gamma}^{\mathrm{solve}}(\mathcal{S}).
\end{equation}
Thus, $\underline{\gamma}^{\mathrm{part}}(\mathcal{S})$ is a lower bound on the
partition gap, whereas $\overline{\gamma}^{\mathrm{solve}}(\mathcal{S})$ is an upper
bound on the subproblem-solution gap.

\section{Computational results}
\label{sec:computational}
It is well known in the literature that ML constraints induce an equivalence relation through their reflexive, symmetric, and transitive closure. Each equivalence class must belong to the same cluster and can therefore be represented by its centroid, weighted by its cardinality. Its within-class sum of squares contributes an additive constant to the objective.
Transferring CL constraints to the corresponding representatives yields an equivalent weighted MSSC problem with only CL constraints, up to this additive constant. Such reductions are used, for example, in \citet{GunsEtAl2016, PiccialliRussoSudoso2022, BaumannHochbaum2025}.
We therefore exclude ML constraints from our experiments, since they can be eliminated in preprocessing. The clustering problem considered here contains only CL constraints ($\mathcal{ML}=\varnothing$).

In this section, we first describe the benchmark instances, the generation of CL
constraint sets, and the experimental setting. We then examine how the number
of anticlusters $T$ affects bound quality and computational performance.
Next, we compare the anticlustering partitions produced by VNS with those obtained by solving
the MISOCP in~\eqref{model:component_ac_miqcp} using a state-of-the-art solver.
Finally, we evaluate the importance of constraint preservation by comparing
the gaps with those obtained from constraint-neglecting partitions.

\subsection{Benchmark datasets and constraints generation}
\label{subsec:instances}
The benchmark comprises four synthetic and six real-world datasets.
The synthetic datasets consist of $K=3$ equally sized Gaussian clusters
in two dimensions. They are denoted by $\texttt{syn-}N\texttt{-}\sigma$,
where $N\in\{1500,3000\}$ is the number of observations and
$\sigma\in\{50,75\}$ controls the common within-cluster dispersion.
The cluster centers are identical across the four datasets, so increasing
$\sigma$ increases cluster overlap. This design allows us to examine
the effects of dataset size and geometric separation.

The real-world datasets are \texttt{Abalone}, \texttt{Rice}, and
\texttt{Wifi} from the UCI Machine Learning Repository;
\texttt{EngyTime} from the Fundamental Clustering Problems Suite;
\texttt{Mallat} from the UCR Time Series Classification Archive;
and \texttt{PR2392} from TSPLIB.
Their configurations $(N,F,K)$, specifying the number of observations,
features, and prescribed clusters, respectively, are
\texttt{Abalone} $(4177,10,3)$,
\texttt{EngyTime} $(4096,2,2)$,
\texttt{Mallat} $(2400,1024,4)$,
\texttt{PR2392} $(2392,2,3)$,
\texttt{Rice} $(3810,7,2)$, and
\texttt{Wifi} $(2000,7,5)$.
Dataset details and download links are provided in Appendix~A.
For each dataset, we retain the partition with the smallest MSSC
objective value among $1{,}000$ random starts of $K$-means, denoted by
$\mathcal P^{(0)}=\{C_1^{(0)},\ldots,C_K^{(0)}\}$.
We then construct a perturbed partition
$\mathcal P^{(1)}=\{C_1^{(1)},\ldots,C_K^{(1)}\}$
by exchanging observations near the boundaries between
clusters with nearby centroids.
This perturbation allows the generated constraints to encode a
partition that differs from the initial $K$-means solution.
Given a density parameter $p\in\{20,40\}$, we generate $\left\lfloor
\frac{p}{100}\min_{j\in[K]}|C_j^{(1)}|
\right\rfloor
\binom{K}{2}$ distinct CL pairs. Each pair connects observations belonging
to different clusters of $\mathcal P^{(1)}$, ensuring that
$\mathcal P^{(1)}$ satisfies all generated constraints and that the
resulting instance is feasible. Pairs are sampled without replacement,
with sampling probabilities favoring observations close to the
decision boundaries of $\mathcal P^{(0)}$ and each observation
may occur in multiple pairs.
For each dataset and each value of $p$, three constraint sets are
generated using random seeds $1$, $2$, and $3$. This yields six
constrained instances per dataset: $24$ synthetic and $36$ real-world
instances, for a total of $60$. 

\subsection{Experimental setting}
\label{subsec:settings}
All experiments are performed on an Apple M4 Max processor with 14 cores and
36 GB of RAM, running macOS 15.5.
For each constrained instance, we use \textsc{PCCC} \citep{BaumannHochbaum2025} to obtain the reference partition $\mathcal P$ and the initial upper bound $UB=f(\mathcal P)$.
For the VNS in Algorithm \ref{alg:vns}, we use $n_{\mathrm{start}}=3$ random starts with a time limit of $60$ seconds each,
$q_{\max}=20$ iterations per start, and maximum shaking strength $\kappa_{\max}=5$. Initialization allows at most $\alpha_\textrm{init}=10{,}000$ attempts per start.
For the CAP cardinality constraints in~\eqref{model:component_ac_cluster_balance}, we set
${\ell}_j(T)=\lfloor0.8\cdot \frac{|C_j|}{T}\rfloor$ and $u_j(T)=\lceil 1.2\cdot \frac{|C_j|}{T}\rceil$.
The constrained MSSC subproblems induced by the anticlusters are handled by \textsc{PC-SOS-SDP} \citep{PiccialliRussoSudoso2022}, a state-of-the-art exact algorithm
that can solve instances with up to 800 observations in presence of ML and CL constraints.
In our experiments, we use 400 observations as a target for the subproblem size.
Bound computation uses a pool of seven workers to process subproblems in parallel, each using two threads for PC-SOS-SDP calls.
For each subproblem, the search is restricted to the root node of the branch-and-bound tree. Subproblem $t$ is declared solved to optimality when the relative gap $\frac{UB_t-LB_t}{UB_t}\le10^{-4}$. All other PC-SOS-SDP parameters are left at their default values.

The choice of $T$ must account for both the subproblem sizes that {PC-SOS-SDP} can handle and the indivisibility of the connected components.
We define $\underline{S} = \max\{100,\lceil1.2 \cdot r_\textrm{max}\rceil\},$ and $\overline{S} = \max\{400,\lceil1.2 \cdot r_\textrm{max}\rceil\},$ where $r_\textrm{max}$ is the size of the largest connected component.
The lower threshold excludes excessively small subproblems, while the upper threshold reflects the practical tractability of {PC-SOS-SDP}. Both thresholds incorporate a $20\%$ margin relative to the largest connected component.
Together, they define the candidate range
$\lceil N/\overline{S} \rceil \le T \le \lfloor N/\underline{S} \rfloor$.
A value of $T$ in this range may still yield an infeasible CAP instance.
We therefore apply the randomized backtracking routine used to initialize
VNS to each candidate and retain $T$ only if the routine constructs a
feasible anticlustering partition.
Let
$\mathcal{T}=\{T_1<T_2<\cdots<T_\tau\}$
denote the ordered set of candidates. We evaluate
$T^{\min}:=T_1$, $T^{\max}:=T_\tau$, and the median
$T^{\mathrm{mid}}:=T_{\lceil \tau/2\rceil}$.

The source code, test instances, generated cannot-link constraint sets, and
reference partitions are available at \url{https://github.com/AnnaLivia/CAP}.

\subsection{Results on synthetic and real instances}
\label{subsec:results}

In this section, we answer three main research questions.
\begin{itemize}
    \item[RQ1.] How does the number of anticlusters affect bound quality and computational effort?
    \item[RQ2.] How does the VNS compare with direct optimization of the MISOCP formulation in terms of partition quality and computational effort?
    \item[RQ3.] How does splitting cannot-link pairs across anticlusters affect bound quality?
\end{itemize}

\begin{figure}
    \centering
    \includegraphics[scale=0.84]{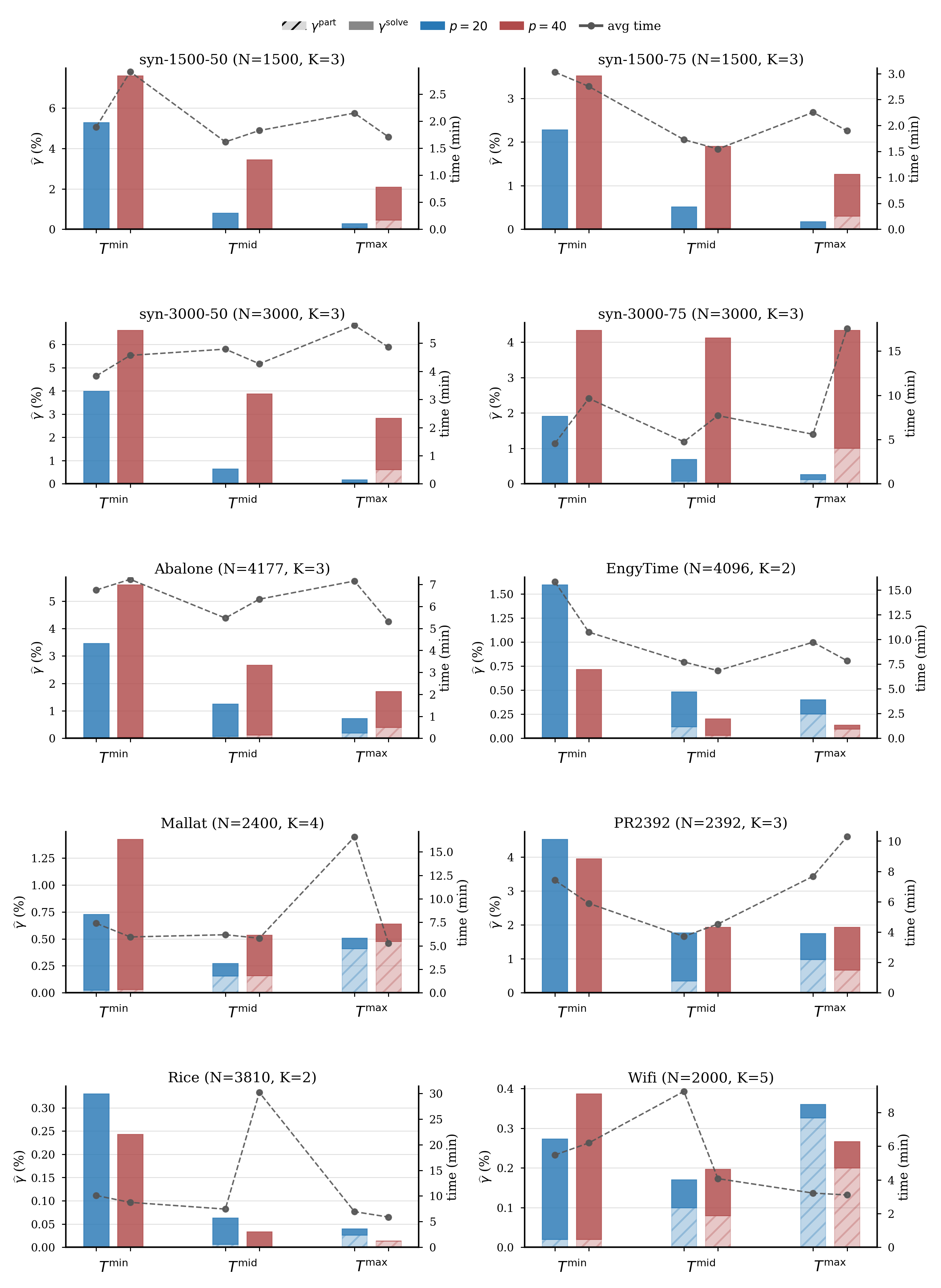}
    \caption{Computable gap bounds and running times for the three anticluster configurations, averaged over the three seeds.}
    \label{fig:t_sensitivity_all}
\end{figure}

\medskip
\paragraph{RQ1}
We compare the three configurations $T^{\mathrm{min}}$, $T^{\mathrm{mid}}$, and $T^{\mathrm{max}}$, which yield progressively smaller constrained MSSC subproblems. Increasing $T$ can make these subproblems easier to solve, but may also weaken the decomposition bound if the smaller anticlusters are less representative of the full dataset. Figure~\ref{fig:t_sensitivity_all} shows the computable gap bounds from~\eqref{eq:observable_gap_bounds} and the corresponding running times, averaged over the three seeds.
For the synthetic instances, moving from $T^{\mathrm{min}}$ to $T^{\mathrm{max}}$ increases the average number of anticlusters from $5$ to $17$ and decreases their average size from $509$ to $173$ observations. The mean gap $\widehat{\gamma}$ decreases from $4.44\%$ to $2.00\%$ and then to $1.43\%$ across the three configurations. For the real-world instances, the average number of anticlusters increases from $8$ to $28$, their average size decreases from $422$ to $118$, and the mean gap decreases from $1.93\%$ to $0.80\%$ and then to $0.71\%$.

These results generally reflect a decrease in $\overline{\gamma}^{\mathrm{solve}}$ accompanied by an increase in $\underline{\gamma}^{\mathrm{part}}$. The former typically outweighs the latter: $T^{\mathrm{max}}$ yields a smaller total gap than $T^{\mathrm{min}}$ in $56$ of the $60$ instances, while $T^{\mathrm{mid}}$ improves on $T^{\mathrm{min}}$ in $59$ instances. For \texttt{Wifi}, \texttt{Mallat}, and \texttt{Rice}, $\underline{\gamma}^{\mathrm{part}}$ accounts for most of the observed gap, but remains small in absolute terms, with maximum values of $0.39\%$, $0.62\%$, and $0.03\%$, respectively.

The effect on running time differs between synthetic and real-world instances. For the latter, the mean lower-bound computation times are $6.28$, $5.75$, and $5.07$ minutes under $T^{\mathrm{min}}$, $T^{\mathrm{mid}}$, and $T^{\mathrm{max}}$, respectively; the corresponding mean total times are $8.14$, $8.13$, and $7.42$ minutes. For the synthetic instances, the mean lower-bound computation times are $3.71$, $2.31$, and $3.72$ minutes, and the mean total times are $4.15$, $3.53$, and $5.21$ minutes. Thus, increasing $T$ generally improves bound quality, but does not necessarily reduce the overall computational cost.

Instances with more CL constraints ($p=40$) have larger mean gaps than those with $p=20$ under all three configurations: $3.44\%$ versus $2.44\%$ for $T^{\mathrm{min}}$, $1.89\%$ versus $0.67\%$ for $T^{\mathrm{mid}}$, and $1.52\%$ versus $0.47\%$ for $T^{\mathrm{max}}$. This pattern is consistent with the role of the largest connected component of the constraint graph. Because each component must be assigned to a single anticluster, a large component imposes a minimum size on that anticluster and, together with the balance constraints, limits the admissible values of $T$. The resulting larger subproblems can be harder to solve, leaving a larger $\overline{\gamma}^{\mathrm{solve}}$, which accounts for most of the observed gap on most instances.

We also consider selecting the best bound among the three tested configurations. For each dataset, constraint density, and seed, let $\mathcal S(T)$ denote the partition obtained with $T$ anticlusters, and select
\begin{equation}\label{eq:tstar}
T^*\in\arg\min_{T\in\{T^{\mathrm{min}},T^{\mathrm{mid}},T^{\mathrm{max}}\}}
\widehat{\gamma}\bigl(\mathcal S(T)\bigr),
\end{equation}
breaking ties in favor of the larger $T$. Since all three configurations yield valid lower bounds and use the same upper bound, this selection returns the strongest of the three computed bounds.

The resulting mean gap is $0.92\%$, compared with $2.94\%$, $1.28\%$, and $0.99\%$ for $T^{\mathrm{min}}$, $T^{\mathrm{mid}}$, and $T^{\mathrm{max}}$, respectively. The mean selected gap is $1.38\%$ for synthetic instances and $0.62\%$ for real-world instances. Although $T^{\mathrm{max}}$ performs best on average, $T^{\mathrm{mid}}$ is selected in two synthetic and seventeen real-world instances. Thus, the largest tested value of $T$ does not always yield the strongest bound.
For each $p$, Table~\ref{tab:best_gap_by_instance} reports the selected gaps.
\begin{table*}[!ht]
    \centering
    \small
    \setlength{\tabcolsep}{4pt}
    \renewcommand{\arraystretch}{1.05}
    \caption{Smallest optimality gaps (\%) over the three tested anticluster configurations for each $p$.}
    \label{tab:best_gap_by_instance}
    \begin{tabular}{lrrr c rrrr c rrrr}
        \toprule
        & & & & & \multicolumn{4}{c}{$p=20$} & & \multicolumn{4}{c}{$p=40$} \\
        Instance & $N$ & $F$ & $K$ & & Seed~1 & Seed~2 & Seed~3 & Avg. & & Seed~1 & Seed~2 & Seed~3 & Avg. \\
        \cline{1-4}\cline{6-9}\cline{11-14}
        \texttt{syn-1500-50} & 1500 & 2 & 3 & & 0.10 & 0.49 & 0.27 & 0.29 & & 3.94 & 1.93 & 0.40 & 2.09 \\
        \texttt{syn-1500-75} & 1500 & 2 & 3 & & 0.22 & 0.19 & 0.12 & 0.18 & & 1.59 & 0.53 & 1.66 & 1.26 \\
        \texttt{syn-3000-50} & 3000 & 2 & 3 & & 0.05 & 0.17 & 0.29 & 0.17 & & 1.06 & 3.24 & 4.18 & 2.83 \\
        \texttt{syn-3000-75} & 3000 & 2 & 3 & & 0.20 & 0.23 & 0.37 & 0.27 & & 3.53 & 4.79 & 3.52 & 3.95 \\
        \texttt{Abalone} & 4177 & 10 & 3 & & 0.76 & 0.74 & 0.66 & 0.72 & & 1.93 & 1.81 & 1.37 & 1.70 \\
        \texttt{EngyTime} & 4096 & 2 & 2 & & 0.36 & 0.29 & 0.39 & 0.35 & & 0.21 & 0.14 & 0.06 & 0.14 \\
        \texttt{Mallat} & 2400 & 1024 & 4 & & 0.26 & 0.26 & 0.30 & 0.27 & & 0.54 & 0.56 & 0.51 & 0.54 \\
        \texttt{PR2392} & 2392 & 2 & 3 & & 1.46 & 1.28 & 1.78 & 1.51 & & 1.31 & 2.14 & 1.88 & 1.78 \\
        \texttt{Rice} & 3810 & 7 & 2 & & 0.04 & 0.04 & 0.04 & 0.04 & & 0.01 & 0.02 & 0.01 & 0.01 \\
        \texttt{Wifi} & 2000 & 7 & 5 & & 0.19 & 0.15 & 0.17 & 0.17 & & 0.20 & 0.23 & 0.14 & 0.19 \\
        \bottomrule
    \end{tabular}
\end{table*}

\medskip
\paragraph{RQ2}
We compare the partition obtained by the VNS heuristic with the incumbent returned by independently solving the MISOCP formulation~\eqref{model:component_ac_miqcp}.
Both methods use the selected number of anticlusters $T^*$ in \eqref{eq:tstar} for each of the $60$ constrained instances. The MISOCP is solved with Gurobi 12 \citep{gurobi} under a time limit of $600$ seconds. Detailed results for each run are reported in Appendix~B.

Of the $60$ runs with Gurobi, $43$ reach the time limit and $17$ terminate because of insufficient memory. These two groups return $37$ and $5$ incumbents, respectively, leaving $18$ runs without a feasible solution. The partition-quality comparison therefore covers the $42$ runs with a valid Gurobi incumbent.
\begin{figure}[!ht]
    \centering
    \includegraphics[scale=0.55]{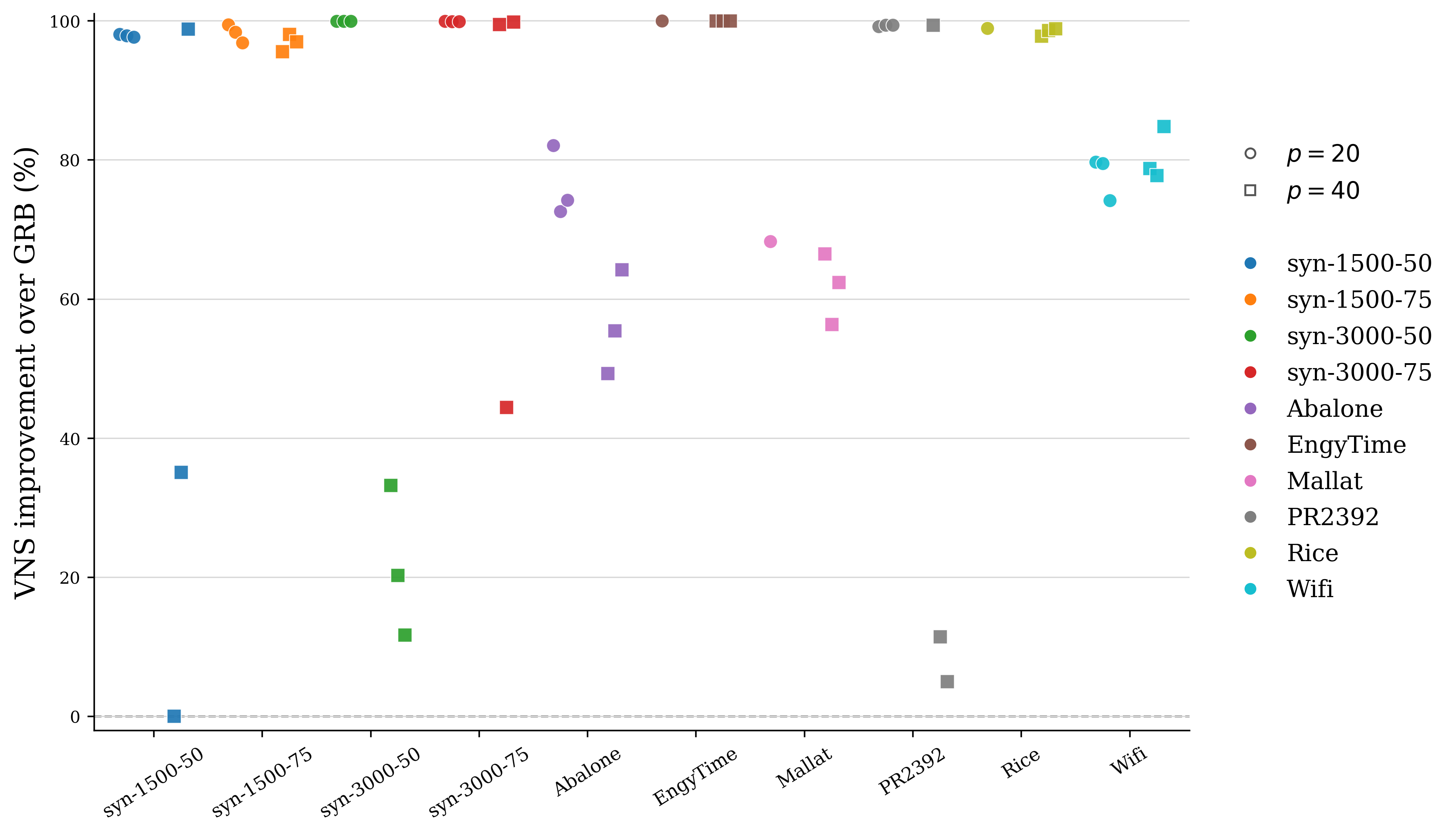}
    \caption{Percentage reduction in CAP objective $D(\mathcal P,\mathcal S)$ achieved by VNS relative to MISOCP solved by Gurobi (GRB).}
    \label{fig:partition_improvement_misocp}
\end{figure}
Figure~\ref{fig:partition_improvement_misocp} reports the percentage reduction in $D(\mathcal P,\mathcal S)$ achieved by VNS with respect to the Gurobi incumbent value (GRB), computed for each matched run as $
\frac{D(\mathcal P,\mathcal S)_{\mathrm{GRB}}
      -D(\mathcal P,\mathcal S)_{\mathrm{VNS}}}
     {D(\mathcal P,\mathcal S)_{\mathrm{GRB}}} \cdot 100$.
Positive values indicate a smaller objective value for the VNS heuristic. VNS outperforms Gurobi in all $42$ matched runs, achieving an $82.2\%$ reduction in the mean objective value. This advantage holds at both constraint densities.
Notably, no Gurobi run certifies optimality. Thus, direct MISOCP optimization struggles to produce high-quality partitions at this scale, whereas VNS consistently returns better partitions whenever a comparison is possible.
In additional runs, we warm-start Gurobi with the VNS solution. Within the same time limit, Gurobi neither improves the initial incumbent nor certifies its optimality. Detailed results for these runs are omitted for brevity.


\medskip
\paragraph{RQ3.}
To assess the effect of preserving CL constraints, we compare CAP with a
constraint-neglecting variant, denoted by CN.
When constructing the anticlusters, CN treats each observation as a singleton component, thereby allowing the endpoints of a CL pair to be assigned to different anticlusters. 
Consequently, a CL constraint is retained in an induced
subproblem only if both endpoints are assigned to the same anticluster;
otherwise, it is omitted from the decomposition. For each dataset,
constraint density, and seed, CN uses the same value $T^*$ selected for CAP.

\begin{figure}[!ht]
    \centering
    \includegraphics[scale=0.5]{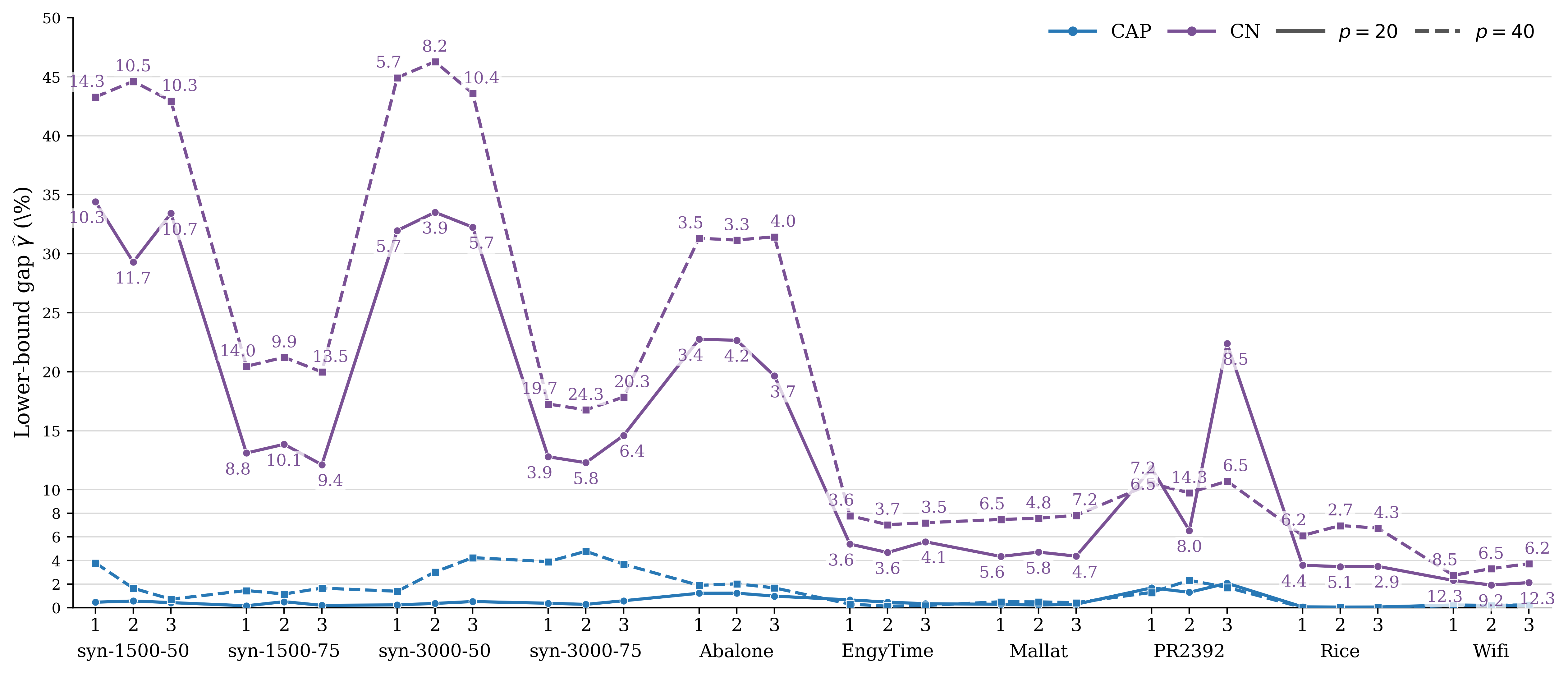}
    \caption{CAP and constraint-neglecting (CN) optimality gaps at the value $T^*$ selected for CAP, averaged over the three seeds. CN labels indicate the percentage of CL constraints retained in the induced subproblems.}
    \label{fig:partition_comparison_gaps}
\end{figure}

Figure~\ref{fig:partition_comparison_gaps} reports the gaps for each dataset and constraint density, averaged over the three seeds. CAP retains all CL constraints by construction, whereas CN retains only $4.64\%$ on average for $p=20$ and $8.04\%$ for $p=40$, or $6.34\%$ overall. CN yields a larger gap in every run. For $p=20$, the average gap increases from $0.40\%$ under CAP to $14.11\%$ under CN; for $p=40$, it increases from $1.45\%$ to $19.69\%$. The deterioration is particularly pronounced for synthetic instances, where the average gap increases from $1.38\%$ to $28.56\%$. For real-world instances, it increases from $0.62\%$ to $9.13\%$.

This comparison also illustrates a limitation of $D(\mathcal P,\mathcal S)$ as a surrogate for bound quality. The objective measures how well the anticlusters reproduce the geometry of the fixed reference partition $\mathcal P$; it does not account for improvements obtained by reassigning observations within the induced subproblems. Splitting CL pairs removes constraints from those subproblems, enlarging their feasible regions and potentially allowing objective values far below those of the restricted reference partitions $\mathcal P|_{S_t}$. Consequently, a small value of $D(\mathcal P,\mathcal S)$ alone does not automatically ensure a strong lower bound. Nevertheless, the results show that preserving CL constraints is important both for bound quality and for the usefulness of $D(\mathcal P,\mathcal S)$ in guiding the decomposition.

\section{Conclusions}
\label{sec:conclusion}
We proposed a decomposition framework for computing lower bounds for large-scale pairwise-constrained MSSC. The framework partitions the observations into anticlusters and combines valid lower bounds from the induced clustering subproblems, providing an optimality certificate for a feasible solution to the original problem. To preserve pairwise information, we introduced the Constrained Anticlustering Problem (CAP), which assigns the connected components of the constraint graph as indivisible units while seeking geometrically representative anticlusters. We formulated CAP as a MISOCP and developed a VNS heuristic to construct high-quality partitions. We also derived computable bounds on the contributions of the decomposition and incomplete subproblem optimization to the resulting optimality gap.
Experiments on synthetic and real-world instances with up to $4{,}177$
observations evaluate the decomposition framework at scales beyond those
reported for direct exact solution with \textsc{PC-SOS-SDP}.
Selecting the strongest bound among the tested anticluster configurations yields a mean optimality gap below $1\%$ at a modest computational cost. VNS produces smaller CAP objective values than direct MISOCP optimization in every run for which the latter returns an incumbent. These results support the use of a heuristic to construct the decomposition, combined with methods for the subproblem to obtain valid bounds.

The experiments also show that geometric representativeness alone is insufficient to obtain strong bounds for constrained clustering. Allowing CL pairs to be split across anticlusters substantially increases the mean optimality gap. Although the resulting bounds remain valid, omitting constraints enlarges the subproblem feasible regions and substantially weakens the bounds. Preserving pairwise information is therefore important for the effectiveness of the decomposition and for the usefulness of the CAP objective as a surrogate for bound quality.

The number of anticlusters controls a trade-off between decomposition quality and subproblem tractability. Larger values generally improve the computed bound by reducing the difficulty of the induced subproblems, but the largest tested value is not always the best. Future work should therefore investigate adaptive strategies that jointly select the number of anticlusters and allocate computational effort across subproblems. Additional improvements may come from alternative anticlustering formulations and more effective exact algorithms for the induced clustering subproblems.


\bibliography{bibfile}

\clearpage

\small
\appendix
\setcounter{table}{0}

\section{Benchmark statistics and CAP results}
\label{appendixA}

Table~\ref{tab:cap_results} summarizes the benchmark statistics and CAP results. The \textit{Dataset} column gives the dataset name and links to its source. The columns $|\mathcal{CL}|$, $R$, and $r_{\max}$ report the number of generated CL pairs, the number of connected components, and the cardinality of the largest component, respectively. The columns $q_{>1}$ and $n_{>1}$ report the number of nonsingleton components and the total number of observations they contain.
The quantities $\mathrm{CL}^{v}_{\text{\tiny UNC}}$ and $\gamma_{\text{\tiny UNC}}$ compare the unconstrained MSSC solution with the corresponding feasible constrained solution for each replicate. Specifically, $\mathrm{CL}^{v}_{\text{\tiny UNC}}$ is the percentage of generated CL constraints violated by the unconstrained solution, while $\gamma_{\text{\tiny UNC}}$ is the relative difference between the objective values of the constrained and unconstrained solutions. Both quantities are expressed as percentages.
For the CAP results, $T$ denotes the number of anticlusters, $|S_t|$ their mean cardinality, and $|\mathcal{CL}_t|$ the mean number of CL constraints retained per anticluster. The quantity $D(\mathcal P,\mathcal S)$ is the aggregate centroid-deviation objective defined in Section~4. The three gap columns report the observable total gap $\widehat{\gamma}$, the lower bound $\underline{\gamma}^{\mathrm{part}}$ on the partition gap, and the upper bound $\overline{\gamma}^{\mathrm{solve}}$ on the subproblem-solution gap. All gaps are expressed as percentages. The final three columns report, in minutes, the running times of the VNS and lower-bound phases and their sum.
For each combination of dataset, constraint level, and replicate, the three rows correspond, in order, to $T^{\min}$, $T^{\mathrm{mid}}$, and $T^{\max}$, with the row attaining the smallest value of $\widehat{\gamma}$ in bold.

\begingroup
\footnotesize
\setlength{\tabcolsep}{1pt}
\renewcommand{\arraystretch}{0.95}
\begin{longtable}{lrr@{\hspace{6pt}}rrrrrrr@{\hspace{6pt}}crrrcrcrrrcrrrc}
\caption{Benchmark statistics and CAP results.}
\label{tab:cap_results} \\
\toprule
\multicolumn{3}{c}{Instance information} & \multicolumn{7}{c}{\hspace{6pt}Benchmark statistics} & \multicolumn{15}{c}{\hspace{6pt}CAP results} \\
\cmidrule(lr){1-3}\cmidrule(lr){4-10}\cmidrule(lr){11-25}
\multicolumn{3}{c}{} & \multicolumn{7}{c}{} & \multicolumn{1}{c}{} & \multicolumn{5}{c}{} & \multicolumn{1}{c}{} & \multicolumn{3}{c}{$\gamma$} & \multicolumn{1}{c}{} & \multicolumn{3}{c}{Time} & \multicolumn{1}{c}{} \\
\cmidrule(lr){18-20}\cmidrule(lr){22-24}
Dataset & $p$ & Rep. & $|\mathcal{CL}|$ & $R$ & $r_{\max}$ & $q_{>1}$ & $n_{>1}$ & $\mathrm{CL}^{v}_{\text{\tiny UNC}}$ & $\gamma_{\text{\tiny UNC}}$ & \multicolumn{1}{c}{} & $T$ & $|S_t|$ & $|\mathcal{CL}_t|$ & \multicolumn{1}{c}{} & $D(\mathcal P,\mathcal S)$ & \multicolumn{1}{c}{} & $\widehat{\gamma}$ & $\underline{\gamma}^{\mathrm{part}}$ & $\overline{\gamma}^{\mathrm{solve}}$ & \multicolumn{1}{c}{} & VNS & LB & TOT & \multicolumn{1}{c}{} \\
\midrule
\endfirsthead
\multicolumn{25}{l}{\tablename~\thetable{} (continued)} \\
\toprule
\multicolumn{3}{c}{Instance information} & \multicolumn{7}{c}{\hspace{6pt}Benchmark statistics} & \multicolumn{15}{c}{\hspace{6pt}CAP results} \\
\cmidrule(lr){1-3}\cmidrule(lr){4-10}\cmidrule(lr){11-25}
\multicolumn{3}{c}{} & \multicolumn{7}{c}{} & \multicolumn{1}{c}{} & \multicolumn{5}{c}{} & \multicolumn{1}{c}{} & \multicolumn{3}{c}{$\gamma$} & \multicolumn{1}{c}{} & \multicolumn{3}{c}{Time} & \multicolumn{1}{c}{} \\
\cmidrule(lr){18-20}\cmidrule(lr){22-24}
Dataset & $p$ & Rep. & $|\mathcal{CL}|$ & $R$ & $r_{\max}$ & $q_{>1}$ & $n_{>1}$ & $\mathrm{CL}^{v}_{\text{\tiny UNC}}$ & $\gamma_{\text{\tiny UNC}}$ & \multicolumn{1}{c}{} & $T$ & $|S_t|$ & $|\mathcal{CL}_t|$ & \multicolumn{1}{c}{} & $D(\mathcal P,\mathcal S)$ & \multicolumn{1}{c}{} & $\widehat{\gamma}$ & $\underline{\gamma}^{\mathrm{part}}$ & $\overline{\gamma}^{\mathrm{solve}}$ & \multicolumn{1}{c}{} & VNS & LB & TOT & \multicolumn{1}{c}{} \\
\midrule
\endhead
\midrule
\multicolumn{25}{r}{Continued on next page} \\
\endfoot
\bottomrule
\endlastfoot
\multirow[c]{18}{*}{\shortstack{\texttt{syn-1500-50}\\\href{https://github.com/AnnaLivia/CAP}{CAP}}} & \multirow[c]{9}{*}{20} & \multirow[c]{3}{*}{1} & \multirow[c]{3}{*}{300} & \multirow[c]{3}{*}{1200} & \multirow[c]{3}{*}{9} & \multirow[c]{3}{*}{166} & \multirow[c]{3}{*}{466} & \multirow[c]{3}{*}{40.00} & \multirow[c]{3}{*}{38.71} &  & 3 & 500 & 100 &  & $<10^{-3}$ &  & 5.51 & 0.00 & 5.51 &  & 0.1 & 1.4 & 1.5 &  \\
 &  &  &  &  &  &  &  &  &  &  & 9 & 167 & 33 &  & $<10^{-3}$ &  & 0.77 & 0.00 & 0.77 &  & 0.7 & 0.9 & 1.6 &  \\
 &  &  &  &  &  &  &  &  &  &  & \textbf{15} & 100 & 20 &  & $<10^{-3}$ &  & 0.10 & 0.02 & 0.08 &  & 1.0 & 1.4 & 2.4 &  \\
\cmidrule(lr){4-10}\cmidrule(lr){11-25}
 &  & \multirow[c]{3}{*}{2} & \multirow[c]{3}{*}{300} & \multirow[c]{3}{*}{1201} & \multirow[c]{3}{*}{20} & \multirow[c]{3}{*}{164} & \multirow[c]{3}{*}{463} & \multirow[c]{3}{*}{35.33} & \multirow[c]{3}{*}{33.89} &  & 3 & 500 & 100 &  & $<10^{-3}$ &  & 4.63 & 0.01 & 4.62 &  & 0.1 & 1.6 & 1.8 &  \\
 &  &  &  &  &  &  &  &  &  &  & 9 & 167 & 33 &  & $<10^{-3}$ &  & 0.65 & 0.03 & 0.62 &  & 0.7 & 0.8 & 1.5 &  \\
 &  &  &  &  &  &  &  &  &  &  & \textbf{15} & 100 & 20 &  & $<10^{-3}$ &  & 0.49 & 0.19 & 0.30 &  & 0.8 & 1.1 & 2.0 &  \\
\cmidrule(lr){4-10}\cmidrule(lr){11-25}
 &  & \multirow[c]{3}{*}{3} & \multirow[c]{3}{*}{300} & \multirow[c]{3}{*}{1200} & \multirow[c]{3}{*}{13} & \multirow[c]{3}{*}{169} & \multirow[c]{3}{*}{469} & \multirow[c]{3}{*}{38.33} & \multirow[c]{3}{*}{38.29} &  & 3 & 500 & 100 &  & $<10^{-3}$ &  & 5.72 & 0.00 & 5.72 &  & 0.1 & 2.3 & 2.4 &  \\
 &  &  &  &  &  &  &  &  &  &  & 9 & 167 & 33 &  & $<10^{-3}$ &  & 1.00 & 0.05 & 0.95 &  & 0.7 & 1.0 & 1.7 &  \\
 &  &  &  &  &  &  &  &  &  &  & \textbf{15} & 100 & 20 &  & $<10^{-3}$ &  & 0.27 & 0.06 & 0.21 &  & 1.0 & 1.1 & 2.1 &  \\
\addlinespace[1pt]
\cmidrule(lr){3-10}\cmidrule(lr){11-25}
 & \multirow[c]{9}{*}{40} & \multirow[c]{3}{*}{1} & \multirow[c]{3}{*}{600} & \multirow[c]{3}{*}{901} & \multirow[c]{3}{*}{170} & \multirow[c]{3}{*}{173} & \multirow[c]{3}{*}{772} & \multirow[c]{3}{*}{35.50} & \multirow[c]{3}{*}{53.03} &  & 3 & 500 & 200 &  & $<10^{-3}$ &  & 8.14 & 0.03 & 8.11 &  & 0.1 & 4.2 & 4.2 &  \\
 &  &  &  &  &  &  &  &  &  &  & 5 & 300 & 120 &  & $<10^{-3}$ &  & 5.90 & 0.07 & 5.83 &  & 0.3 & 2.2 & 2.4 &  \\
 &  &  &  &  &  &  &  &  &  &  & \textbf{8} & 188 & 75 &  & 7.901 &  & 3.94 & 0.79 & 3.15 &  & 0.3 & 1.1 & 1.5 &  \\
\cmidrule(lr){4-10}\cmidrule(lr){11-25}
 &  & \multirow[c]{3}{*}{2} & \multirow[c]{3}{*}{600} & \multirow[c]{3}{*}{903} & \multirow[c]{3}{*}{117} & \multirow[c]{3}{*}{169} & \multirow[c]{3}{*}{766} & \multirow[c]{3}{*}{35.83} & \multirow[c]{3}{*}{50.88} &  & 3 & 500 & 200 &  & $<10^{-3}$ &  & 8.15 & 0.00 & 8.15 &  & 0.1 & 1.7 & 1.8 &  \\
 &  &  &  &  &  &  &  &  &  &  & 7 & 214 & 86 &  & $<10^{-3}$ &  & 2.99 & 0.00 & 2.99 &  & 0.3 & 1.4 & 1.7 &  \\
 &  &  &  &  &  &  &  &  &  &  & \textbf{11} & 136 & 55 &  & 7.873 &  & 1.93 & 0.55 & 1.38 &  & 0.4 & 1.2 & 1.6 &  \\
\cmidrule(lr){4-10}\cmidrule(lr){11-25}
 &  & \multirow[c]{3}{*}{3} & \multirow[c]{3}{*}{600} & \multirow[c]{3}{*}{900} & \multirow[c]{3}{*}{87} & \multirow[c]{3}{*}{183} & \multirow[c]{3}{*}{783} & \multirow[c]{3}{*}{33.67} & \multirow[c]{3}{*}{48.89} &  & 3 & 500 & 200 &  & $<10^{-3}$ &  & 6.48 & 0.00 & 6.48 &  & 0.1 & 2.6 & 2.7 &  \\
 &  &  &  &  &  &  &  &  &  &  & 9 & 167 & 67 &  & $<10^{-3}$ &  & 1.43 & 0.00 & 1.43 &  & 0.4 & 1.0 & 1.4 &  \\
 &  &  &  &  &  &  &  &  &  &  & \textbf{15} & 100 & 40 &  & 0.474 &  & 0.40 & 0.03 & 0.37 &  & 0.5 & 1.5 & 2.0 &  \\
\addlinespace[1pt]
\pagebreak
\multirow[c]{18}{*}{\shortstack{\texttt{syn-1500-75}\\\href{https://github.com/AnnaLivia/CAP}{CAP}}} & \multirow[c]{9}{*}{20} & \multirow[c]{3}{*}{1} & \multirow[c]{3}{*}{297} & \multirow[c]{3}{*}{1203} & \multirow[c]{3}{*}{19} & \multirow[c]{3}{*}{138} & \multirow[c]{3}{*}{435} & \multirow[c]{3}{*}{39.39} & \multirow[c]{3}{*}{14.11} &  & 3 & 500 & 99 &  & $<10^{-3}$ &  & 2.14 & 0.00 & 2.14 &  & 0.2 & 2.9 & 3.1 &  \\
 &  &  &  &  &  &  &  &  &  &  & 9 & 167 & 33 &  & $<10^{-3}$ &  & 0.38 & 0.01 & 0.37 &  & 0.7 & 0.9 & 1.6 &  \\
 &  &  &  &  &  &  &  &  &  &  & \textbf{15} & 100 & 20 &  & $<10^{-3}$ &  & 0.22 & 0.04 & 0.18 &  & 0.9 & 1.4 & 2.3 &  \\
\cmidrule(lr){4-10}\cmidrule(lr){11-25}
 &  & \multirow[c]{3}{*}{2} & \multirow[c]{3}{*}{297} & \multirow[c]{3}{*}{1203} & \multirow[c]{3}{*}{16} & \multirow[c]{3}{*}{155} & \multirow[c]{3}{*}{452} & \multirow[c]{3}{*}{38.38} & \multirow[c]{3}{*}{14.80} &  & 3 & 500 & 99 &  & $<10^{-3}$ &  & 2.64 & 0.00 & 2.64 &  & 0.1 & 3.3 & 3.5 &  \\
 &  &  &  &  &  &  &  &  &  &  & 9 & 167 & 33 &  & $<10^{-3}$ &  & 0.91 & 0.00 & 0.91 &  & 0.7 & 1.1 & 1.8 &  \\
 &  &  &  &  &  &  &  &  &  &  & \textbf{15} & 100 & 20 &  & $<10^{-3}$ &  & 0.19 & 0.01 & 0.18 &  & 0.9 & 1.3 & 2.2 &  \\
\cmidrule(lr){4-10}\cmidrule(lr){11-25}
 &  & \multirow[c]{3}{*}{3} & \multirow[c]{3}{*}{297} & \multirow[c]{3}{*}{1203} & \multirow[c]{3}{*}{20} & \multirow[c]{3}{*}{135} & \multirow[c]{3}{*}{432} & \multirow[c]{3}{*}{35.35} & \multirow[c]{3}{*}{13.40} &  & 3 & 500 & 99 &  & $<10^{-3}$ &  & 2.05 & 0.00 & 2.05 &  & 0.1 & 2.4 & 2.6 &  \\
 &  &  &  &  &  &  &  &  &  &  & 9 & 167 & 33 &  & $<10^{-3}$ &  & 0.25 & 0.00 & 0.25 &  & 0.7 & 1.0 & 1.8 &  \\
 &  &  &  &  &  &  &  &  &  &  & \textbf{15} & 100 & 20 &  & $<10^{-3}$ &  & 0.12 & 0.00 & 0.12 &  & 1.0 & 1.3 & 2.3 &  \\
\addlinespace[1pt]
\cmidrule(lr){3-10}\cmidrule(lr){11-25}
 & \multirow[c]{9}{*}{40} & \multirow[c]{3}{*}{1} & \multirow[c]{3}{*}{594} & \multirow[c]{3}{*}{908} & \multirow[c]{3}{*}{140} & \multirow[c]{3}{*}{151} & \multirow[c]{3}{*}{743} & \multirow[c]{3}{*}{38.89} & \multirow[c]{3}{*}{24.02} &  & 3 & 500 & 198 &  & $<10^{-3}$ &  & 3.83 & 0.00 & 3.83 &  & 0.1 & 1.7 & 1.8 &  \\
 &  &  &  &  &  &  &  &  &  &  & 6 & 250 & 99 &  & $<10^{-3}$ &  & 2.19 & 0.00 & 2.19 &  & 0.3 & 1.3 & 1.6 &  \\
 &  &  &  &  &  &  &  &  &  &  & \textbf{10} & 150 & 59 &  & 6.952 &  & 1.59 & 0.55 & 1.04 &  & 0.5 & 1.0 & 1.5 &  \\
\cmidrule(lr){4-10}\cmidrule(lr){11-25}
 &  & \multirow[c]{3}{*}{2} & \multirow[c]{3}{*}{594} & \multirow[c]{3}{*}{907} & \multirow[c]{3}{*}{78} & \multirow[c]{3}{*}{154} & \multirow[c]{3}{*}{747} & \multirow[c]{3}{*}{36.03} & \multirow[c]{3}{*}{23.08} &  & 3 & 500 & 198 &  & $<10^{-3}$ &  & 3.69 & 0.00 & 3.69 &  & 0.1 & 2.7 & 2.8 &  \\
 &  &  &  &  &  &  &  &  &  &  & 9 & 167 & 66 &  & $<10^{-3}$ &  & 1.16 & 0.06 & 1.10 &  & 0.4 & 1.0 & 1.4 &  \\
 &  &  &  &  &  &  &  &  &  &  & \textbf{15} & 100 & 40 &  & 0.001 &  & 0.53 & 0.16 & 0.37 &  & 0.5 & 2.1 & 2.6 &  \\
\cmidrule(lr){4-10}\cmidrule(lr){11-25}
 &  & \multirow[c]{3}{*}{3} & \multirow[c]{3}{*}{594} & \multirow[c]{3}{*}{909} & \multirow[c]{3}{*}{168} & \multirow[c]{3}{*}{144} & \multirow[c]{3}{*}{735} & \multirow[c]{3}{*}{37.21} & \multirow[c]{3}{*}{23.50} &  & 3 & 500 & 198 &  & $<10^{-3}$ &  & 3.02 & 0.00 & 3.02 &  & 0.1 & 3.7 & 3.7 &  \\
 &  &  &  &  &  &  &  &  &  &  & 5 & 300 & 119 &  & $<10^{-3}$ &  & 2.35 & 0.04 & 2.31 &  & 0.3 & 1.3 & 1.6 &  \\
 &  &  &  &  &  &  &  &  &  &  & \textbf{8} & 188 & 74 &  & 2.749 &  & 1.66 & 0.20 & 1.46 &  & 0.4 & 1.2 & 1.6 &  \\
\addlinespace[1pt]
\midrule
\multirow[c]{18}{*}{\shortstack{\texttt{syn-3000-50}\\\href{https://github.com/AnnaLivia/CAP}{CAP}}} & \multirow[c]{9}{*}{20} & \multirow[c]{3}{*}{1} & \multirow[c]{3}{*}{597} & \multirow[c]{3}{*}{2403} & \multirow[c]{3}{*}{17} & \multirow[c]{3}{*}{354} & \multirow[c]{3}{*}{951} & \multirow[c]{3}{*}{33.84} & \multirow[c]{3}{*}{34.47} &  & 7 & 429 & 85 &  & $<10^{-3}$ &  & 3.86 & 0.00 & 3.86 &  & 1.0 & 2.5 & 3.5 &  \\
 &  &  &  &  &  &  &  &  &  &  & 18 & 167 & 33 &  & $<10^{-3}$ &  & 0.43 & 0.01 & 0.42 &  & 3.0 & 1.6 & 4.7 &  \\
 &  &  &  &  &  &  &  &  &  &  & \textbf{30} & 100 & 20 &  & $<10^{-3}$ &  & 0.05 & 0.01 & 0.04 &  & 3.2 & 2.4 & 5.6 &  \\
\cmidrule(lr){4-10}\cmidrule(lr){11-25}
 &  & \multirow[c]{3}{*}{2} & \multirow[c]{3}{*}{597} & \multirow[c]{3}{*}{2403} & \multirow[c]{3}{*}{15} & \multirow[c]{3}{*}{324} & \multirow[c]{3}{*}{921} & \multirow[c]{3}{*}{33.00} & \multirow[c]{3}{*}{34.88} &  & 7 & 429 & 85 &  & $<10^{-3}$ &  & 4.12 & 0.00 & 4.12 &  & 1.1 & 2.8 & 3.9 &  \\
 &  &  &  &  &  &  &  &  &  &  & 18 & 167 & 33 &  & $<10^{-3}$ &  & 0.68 & 0.01 & 0.67 &  & 3.0 & 1.9 & 4.9 &  \\
 &  &  &  &  &  &  &  &  &  &  & \textbf{30} & 100 & 20 &  & $<10^{-3}$ &  & 0.17 & 0.07 & 0.10 &  & 3.1 & 2.4 & 5.4 &  \\
\cmidrule(lr){4-10}\cmidrule(lr){11-25}
 &  & \multirow[c]{3}{*}{3} & \multirow[c]{3}{*}{597} & \multirow[c]{3}{*}{2403} & \multirow[c]{3}{*}{14} & \multirow[c]{3}{*}{335} & \multirow[c]{3}{*}{932} & \multirow[c]{3}{*}{32.50} & \multirow[c]{3}{*}{33.14} &  & 7 & 429 & 85 &  & $<10^{-3}$ &  & 4.00 & 0.00 & 4.00 &  & 1.1 & 3.0 & 4.1 &  \\
 &  &  &  &  &  &  &  &  &  &  & 18 & 167 & 33 &  & $<10^{-3}$ &  & 0.80 & 0.03 & 0.77 &  & 3.0 & 1.8 & 4.8 &  \\
 &  &  &  &  &  &  &  &  &  &  & \textbf{30} & 100 & 20 &  & $<10^{-3}$ &  & 0.29 & 0.11 & 0.18 &  & 3.1 & 2.7 & 5.8 &  \\
\addlinespace[1pt]
\cmidrule(lr){3-10}\cmidrule(lr){11-25}
 & \multirow[c]{9}{*}{40} & \multirow[c]{3}{*}{1} & \multirow[c]{3}{*}{1197} & \multirow[c]{3}{*}{1803} & \multirow[c]{3}{*}{129} & \multirow[c]{3}{*}{328} & \multirow[c]{3}{*}{1525} & \multirow[c]{3}{*}{33.58} & \multirow[c]{3}{*}{48.10} &  & 7 & 429 & 171 &  & $<10^{-3}$ &  & 7.14 & 0.00 & 7.14 &  & 0.6 & 3.1 & 3.7 &  \\
 &  &  &  &  &  &  &  &  &  &  & 14 & 214 & 86 &  & $<10^{-3}$ &  & 2.68 & 0.01 & 2.67 &  & 1.0 & 2.6 & 3.6 &  \\
 &  &  &  &  &  &  &  &  &  &  & \textbf{21} & 143 & 57 &  & 3.957 &  & 1.06 & 0.17 & 0.89 &  & 3.0 & 2.0 & 5.0 &  \\
\cmidrule(lr){4-10}\cmidrule(lr){11-25}
 &  & \multirow[c]{3}{*}{2} & \multirow[c]{3}{*}{1197} & \multirow[c]{3}{*}{1804} & \multirow[c]{3}{*}{159} & \multirow[c]{3}{*}{336} & \multirow[c]{3}{*}{1532} & \multirow[c]{3}{*}{35.92} & \multirow[c]{3}{*}{49.01} &  & 7 & 429 & 171 &  & $<10^{-3}$ &  & 6.57 & 0.00 & 6.57 &  & 0.6 & 5.0 & 5.6 &  \\
 &  &  &  &  &  &  &  &  &  &  & 12 & 250 & 100 &  & $<10^{-3}$ &  & 4.31 & 0.02 & 4.29 &  & 1.9 & 3.0 & 4.9 &  \\
 &  &  &  &  &  &  &  &  &  &  & \textbf{18} & 167 & 67 &  & 22.281 &  & 3.24 & 0.83 & 2.41 &  & 2.5 & 2.1 & 4.6 &  \\
\cmidrule(lr){4-10}\cmidrule(lr){11-25}
 &  & \multirow[c]{3}{*}{3} & \multirow[c]{3}{*}{1197} & \multirow[c]{3}{*}{1804} & \multirow[c]{3}{*}{231} & \multirow[c]{3}{*}{322} & \multirow[c]{3}{*}{1518} & \multirow[c]{3}{*}{35.25} & \multirow[c]{3}{*}{48.82} &  & 7 & 429 & 171 &  & $<10^{-3}$ &  & 6.12 & 0.01 & 6.11 &  & 0.5 & 3.9 & 4.4 &  \\
 &  &  &  &  &  &  &  &  &  &  & 9 & 333 & 133 &  & 0.062 &  & 4.66 & 0.03 & 4.63 &  & 0.7 & 3.6 & 4.3 &  \\
 &  &  &  &  &  &  &  &  &  &  & \textbf{12} & 250 & 100 &  & 18.300 &  & 4.18 & 0.85 & 3.33 &  & 2.0 & 3.0 & 5.0 &  \\
\addlinespace[1pt]
\pagebreak
\multirow[c]{18}{*}{\shortstack{\texttt{syn-3000-75}\\\href{https://github.com/AnnaLivia/CAP}{CAP}}} & \multirow[c]{9}{*}{20} & \multirow[c]{3}{*}{1} & \multirow[c]{3}{*}{591} & \multirow[c]{3}{*}{2409} & \multirow[c]{3}{*}{18} & \multirow[c]{3}{*}{317} & \multirow[c]{3}{*}{908} & \multirow[c]{3}{*}{34.69} & \multirow[c]{3}{*}{12.67} &  & 7 & 429 & 84 &  & $<10^{-3}$ &  & 1.78 & 0.01 & 1.77 &  & 1.0 & 3.3 & 4.3 &  \\
 &  &  &  &  &  &  &  &  &  &  & 18 & 167 & 33 &  & $<10^{-3}$ &  & 0.61 & 0.04 & 0.57 &  & 3.0 & 1.7 & 4.7 &  \\
 &  &  &  &  &  &  &  &  &  &  & \textbf{30} & 100 & 20 &  & 0.001 &  & 0.20 & 0.11 & 0.09 &  & 2.3 & 2.5 & 4.8 &  \\
\cmidrule(lr){4-10}\cmidrule(lr){11-25}
 &  & \multirow[c]{3}{*}{2} & \multirow[c]{3}{*}{591} & \multirow[c]{3}{*}{2410} & \multirow[c]{3}{*}{16} & \multirow[c]{3}{*}{321} & \multirow[c]{3}{*}{911} & \multirow[c]{3}{*}{36.38} & \multirow[c]{3}{*}{12.62} &  & 7 & 429 & 84 &  & $<10^{-3}$ &  & 1.61 & 0.00 & 1.61 &  & 1.1 & 3.2 & 4.3 &  \\
 &  &  &  &  &  &  &  &  &  &  & 18 & 167 & 33 &  & $<10^{-3}$ &  & 0.54 & 0.06 & 0.48 &  & 3.0 & 1.7 & 4.7 &  \\
 &  &  &  &  &  &  &  &  &  &  & \textbf{30} & 100 & 20 &  & 0.001 &  & 0.23 & 0.14 & 0.09 &  & 3.0 & 2.8 & 5.8 &  \\
\cmidrule(lr){4-10}\cmidrule(lr){11-25}
 &  & \multirow[c]{3}{*}{3} & \multirow[c]{3}{*}{591} & \multirow[c]{3}{*}{2409} & \multirow[c]{3}{*}{15} & \multirow[c]{3}{*}{329} & \multirow[c]{3}{*}{920} & \multirow[c]{3}{*}{38.07} & \multirow[c]{3}{*}{15.00} &  & 7 & 429 & 84 &  & $<10^{-3}$ &  & 2.33 & 0.01 & 2.32 &  & 1.0 & 4.0 & 5.0 &  \\
 &  &  &  &  &  &  &  &  &  &  & 18 & 167 & 33 &  & $<10^{-3}$ &  & 0.91 & 0.12 & 0.79 &  & 3.0 & 1.8 & 4.8 &  \\
 &  &  &  &  &  &  &  &  &  &  & \textbf{30} & 100 & 20 &  & $<10^{-3}$ &  & 0.37 & 0.13 & 0.24 &  & 3.0 & 3.1 & 6.2 &  \\
\addlinespace[1pt]
\cmidrule(lr){3-10}\cmidrule(lr){11-25}
 & \multirow[c]{9}{*}{40} & \multirow[c]{3}{*}{1} & \multirow[c]{3}{*}{1182} & \multirow[c]{3}{*}{1822} & \multirow[c]{3}{*}{418} & \multirow[c]{3}{*}{265} & \multirow[c]{3}{*}{1443} & \multirow[c]{3}{*}{37.06} & \multirow[c]{3}{*}{21.61} &  & 5 & 600 & 236 &  & $<10^{-3}$ &  & 3.90 & 0.00 & 3.90 &  & 0.4 & 5.9 & 6.3 &  \\
 &  &  &  &  &  &  &  &  &  &  & \textbf{6} & 500 & 197 &  & 3.700 &  & 3.53 & 0.10 & 3.43 &  & 0.6 & 4.4 & 5.0 &  \\
 &  &  &  &  &  &  &  &  &  &  & 7 & 429 & 169 &  & 46.735 &  & 4.35 & 1.35 & 3.00 &  & 0.6 & 38.4 & 39.0 &  \\
\cmidrule(lr){4-10}\cmidrule(lr){11-25}
 &  & \multirow[c]{3}{*}{2} & \multirow[c]{3}{*}{1182} & \multirow[c]{3}{*}{1826} & \multirow[c]{3}{*}{572} & \multirow[c]{3}{*}{272} & \multirow[c]{3}{*}{1446} & \multirow[c]{3}{*}{37.14} & \multirow[c]{3}{*}{21.58} &  & 3 & 1000 & 394 &  & $<10^{-3}$ &  & 5.15 & 0.00 & 5.15 &  & 0.2 & 10.8 & 11.1 &  \\
 &  &  &  &  &  &  &  &  &  &  & \textbf{4} & 750 & 296 &  & 0.010 &  & 4.79 & 0.00 & 4.79 &  & 0.4 & 8.9 & 9.3 &  \\
 &  &  &  &  &  &  &  &  &  &  & 5 & 600 & 236 &  & 54.913 &  & 5.12 & 1.38 & 3.74 &  & 0.4 & 6.6 & 7.0 &  \\
\cmidrule(lr){4-10}\cmidrule(lr){11-25}
 &  & \multirow[c]{3}{*}{3} & \multirow[c]{3}{*}{1182} & \multirow[c]{3}{*}{1819} & \multirow[c]{3}{*}{433} & \multirow[c]{3}{*}{287} & \multirow[c]{3}{*}{1468} & \multirow[c]{3}{*}{36.13} & \multirow[c]{3}{*}{22.05} &  & 4 & 750 & 296 &  & $<10^{-3}$ &  & 3.94 & 0.00 & 3.94 &  & 0.3 & 11.2 & 11.5 &  \\
 &  &  &  &  &  &  &  &  &  &  & 5 & 600 & 236 &  & $<10^{-3}$ &  & 4.03 & 0.00 & 4.03 &  & 0.4 & 8.3 & 8.8 &  \\
 &  &  &  &  &  &  &  &  &  &  & \textbf{6} & 500 & 197 &  & 11.432 &  & 3.52 & 0.29 & 3.23 &  & 1.2 & 5.4 & 6.7 &  \\
\addlinespace[1pt]
\midrule
\multirow[c]{18}{*}{\shortstack{\texttt{Abalone}\\\href{https://archive.ics.uci.edu/dataset/1/abalone}{UCI}}} & \multirow[c]{9}{*}{20} & \multirow[c]{3}{*}{1} & \multirow[c]{3}{*}{783} & \multirow[c]{3}{*}{3394} & \multirow[c]{3}{*}{12} & \multirow[c]{3}{*}{512} & \multirow[c]{3}{*}{1295} & \multirow[c]{3}{*}{34.10} & \multirow[c]{3}{*}{23.56} &  & 10 & 418 & 78 &  & 0.049 &  & 3.69 & 0.03 & 3.66 &  & 3.0 & 4.3 & 7.3 &  \\
 &  &  &  &  &  &  &  &  &  &  & 24 & 174 & 33 &  & 0.428 &  & 1.31 & 0.09 & 1.22 &  & 3.0 & 2.5 & 5.5 &  \\
 &  &  &  &  &  &  &  &  &  &  & \textbf{39} & 107 & 20 &  & 1.101 &  & 0.76 & 0.20 & 0.56 &  & 2.8 & 4.1 & 6.9 &  \\
\cmidrule(lr){4-10}\cmidrule(lr){11-25}
 &  & \multirow[c]{3}{*}{2} & \multirow[c]{3}{*}{783} & \multirow[c]{3}{*}{3394} & \multirow[c]{3}{*}{11} & \multirow[c]{3}{*}{504} & \multirow[c]{3}{*}{1287} & \multirow[c]{3}{*}{35.25} & \multirow[c]{3}{*}{23.76} &  & 10 & 418 & 78 &  & 0.058 &  & 3.77 & 0.04 & 3.73 &  & 3.0 & 3.3 & 6.3 &  \\
 &  &  &  &  &  &  &  &  &  &  & 25 & 167 & 31 &  & 0.504 &  & 1.36 & 0.08 & 1.28 &  & 3.0 & 2.4 & 5.4 &  \\
 &  &  &  &  &  &  &  &  &  &  & \textbf{40} & 104 & 20 &  & 1.101 &  & 0.74 & 0.20 & 0.54 &  & 2.8 & 4.8 & 7.6 &  \\
\cmidrule(lr){4-10}\cmidrule(lr){11-25}
 &  & \multirow[c]{3}{*}{3} & \multirow[c]{3}{*}{783} & \multirow[c]{3}{*}{3394} & \multirow[c]{3}{*}{10} & \multirow[c]{3}{*}{452} & \multirow[c]{3}{*}{1235} & \multirow[c]{3}{*}{30.40} & \multirow[c]{3}{*}{20.43} &  & 10 & 418 & 78 &  & 0.055 &  & 2.91 & 0.01 & 2.90 &  & 3.0 & 3.7 & 6.7 &  \\
 &  &  &  &  &  &  &  &  &  &  & 24 & 174 & 33 &  & 0.350 &  & 1.09 & 0.06 & 1.03 &  & 3.0 & 2.5 & 5.5 &  \\
 &  &  &  &  &  &  &  &  &  &  & \textbf{39} & 107 & 20 &  & 1.115 &  & 0.66 & 0.17 & 0.49 &  & 3.4 & 3.6 & 6.9 &  \\
\addlinespace[1pt]
\cmidrule(lr){3-10}\cmidrule(lr){11-25}
 & \multirow[c]{9}{*}{40} & \multirow[c]{3}{*}{1} & \multirow[c]{3}{*}{1569} & \multirow[c]{3}{*}{2609} & \multirow[c]{3}{*}{127} & \multirow[c]{3}{*}{538} & \multirow[c]{3}{*}{2106} & \multirow[c]{3}{*}{33.01} & \multirow[c]{3}{*}{32.63} &  & 10 & 418 & 157 &  & 0.050 &  & 5.48 & 0.01 & 5.47 &  & 2.9 & 3.8 & 6.6 &  \\
 &  &  &  &  &  &  &  &  &  &  & 20 & 209 & 78 &  & 1.275 &  & 2.62 & 0.12 & 2.50 &  & 3.0 & 3.3 & 6.3 &  \\
 &  &  &  &  &  &  &  &  &  &  & \textbf{31} & 135 & 51 &  & 5.804 &  & 1.93 & 0.49 & 1.44 &  & 1.5 & 2.6 & 4.1 &  \\
\cmidrule(lr){4-10}\cmidrule(lr){11-25}
 &  & \multirow[c]{3}{*}{2} & \multirow[c]{3}{*}{1569} & \multirow[c]{3}{*}{2608} & \multirow[c]{3}{*}{139} & \multirow[c]{3}{*}{544} & \multirow[c]{3}{*}{2113} & \multirow[c]{3}{*}{31.55} & \multirow[c]{3}{*}{32.04} &  & 10 & 418 & 157 &  & 0.046 &  & 5.60 & 0.07 & 5.53 &  & 3.0 & 4.8 & 7.8 &  \\
 &  &  &  &  &  &  &  &  &  &  & 19 & 220 & 83 &  & 0.496 &  & 2.83 & 0.13 & 2.70 &  & 3.0 & 3.6 & 6.6 &  \\
 &  &  &  &  &  &  &  &  &  &  & \textbf{28} & 149 & 56 &  & 2.352 &  & 1.81 & 0.29 & 1.52 &  & 3.4 & 2.6 & 6.0 &  \\
\cmidrule(lr){4-10}\cmidrule(lr){11-25}
 &  & \multirow[c]{3}{*}{3} & \multirow[c]{3}{*}{1569} & \multirow[c]{3}{*}{2609} & \multirow[c]{3}{*}{100} & \multirow[c]{3}{*}{547} & \multirow[c]{3}{*}{2115} & \multirow[c]{3}{*}{32.76} & \multirow[c]{3}{*}{33.47} &  & 10 & 418 & 157 &  & 0.054 &  & 5.69 & 0.03 & 5.66 &  & 2.9 & 4.5 & 7.3 &  \\
 &  &  &  &  &  &  &  &  &  &  & 22 & 190 & 71 &  & 0.718 &  & 2.57 & 0.13 & 2.44 &  & 3.0 & 3.1 & 6.1 &  \\
 &  &  &  &  &  &  &  &  &  &  & \textbf{34} & 123 & 46 &  & 3.722 &  & 1.37 & 0.43 & 0.94 &  & 2.5 & 3.3 & 5.8 &  \\
\addlinespace[1pt]
\pagebreak
\multirow[c]{18}{*}{\shortstack{\texttt{EngyTime}\\\href{https://data.mendeley.com/datasets/vsxvgc4rwy/1}{Mendeley}}} & \multirow[c]{9}{*}{20} & \multirow[c]{3}{*}{1} & \multirow[c]{3}{*}{388} & \multirow[c]{3}{*}{3708} & \multirow[c]{3}{*}{7} & \multirow[c]{3}{*}{293} & \multirow[c]{3}{*}{681} & \multirow[c]{3}{*}{51.55} & \multirow[c]{3}{*}{3.58} &  & 10 & 410 & 39 &  & $<10^{-3}$ &  & 1.84 & 0.01 & 1.83 &  & 1.5 & 12.7 & 14.2 &  \\
 &  &  &  &  &  &  &  &  &  &  & 25 & 164 & 16 &  & $<10^{-3}$ &  & 0.53 & 0.15 & 0.38 &  & 2.8 & 4.8 & 7.5 &  \\
 &  &  &  &  &  &  &  &  &  &  & \textbf{40} & 102 & 10 &  & $<10^{-3}$ &  & 0.36 & 0.25 & 0.11 &  & 3.0 & 7.3 & 10.3 &  \\
\cmidrule(lr){4-10}\cmidrule(lr){11-25}
 &  & \multirow[c]{3}{*}{2} & \multirow[c]{3}{*}{388} & \multirow[c]{3}{*}{3708} & \multirow[c]{3}{*}{11} & \multirow[c]{3}{*}{284} & \multirow[c]{3}{*}{672} & \multirow[c]{3}{*}{57.47} & \multirow[c]{3}{*}{3.54} &  & 10 & 410 & 39 &  & $<10^{-3}$ &  & 1.38 & 0.00 & 1.38 &  & 1.6 & 16.7 & 18.3 &  \\
 &  &  &  &  &  &  &  &  &  &  & 25 & 164 & 16 &  & $<10^{-3}$ &  & 0.53 & 0.14 & 0.39 &  & 3.0 & 4.5 & 7.5 &  \\
 &  &  &  &  &  &  &  &  &  &  & \textbf{40} & 102 & 10 &  & $<10^{-3}$ &  & 0.29 & 0.09 & 0.20 &  & 3.0 & 6.4 & 9.4 &  \\
\cmidrule(lr){4-10}\cmidrule(lr){11-25}
 &  & \multirow[c]{3}{*}{3} & \multirow[c]{3}{*}{388} & \multirow[c]{3}{*}{3708} & \multirow[c]{3}{*}{10} & \multirow[c]{3}{*}{288} & \multirow[c]{3}{*}{676} & \multirow[c]{3}{*}{56.19} & \multirow[c]{3}{*}{3.57} &  & 10 & 410 & 39 &  & $<10^{-3}$ &  & 1.56 & 0.00 & 1.56 &  & 1.5 & 13.5 & 15.1 &  \\
 &  &  &  &  &  &  &  &  &  &  & \textbf{25} & 164 & 16 &  & $<10^{-3}$ &  & 0.39 & 0.07 & 0.32 &  & 3.0 & 5.1 & 8.1 &  \\
 &  &  &  &  &  &  &  &  &  &  & 40 & 102 & 10 &  & $<10^{-3}$ &  & 0.55 & 0.42 & 0.13 &  & 3.0 & 6.5 & 9.5 &  \\
\addlinespace[1pt]
\cmidrule(lr){3-10}\cmidrule(lr){11-25}
 & \multirow[c]{9}{*}{40} & \multirow[c]{3}{*}{1} & \multirow[c]{3}{*}{776} & \multirow[c]{3}{*}{3320} & \multirow[c]{3}{*}{23} & \multirow[c]{3}{*}{405} & \multirow[c]{3}{*}{1181} & \multirow[c]{3}{*}{54.25} & \multirow[c]{3}{*}{6.69} &  & 10 & 410 & 78 &  & $<10^{-3}$ &  & 0.74 & 0.00 & 0.74 &  & 1.3 & 8.6 & 9.9 &  \\
 &  &  &  &  &  &  &  &  &  &  & 25 & 164 & 31 &  & $<10^{-3}$ &  & 0.26 & 0.03 & 0.23 &  & 2.3 & 4.0 & 6.3 &  \\
 &  &  &  &  &  &  &  &  &  &  & \textbf{40} & 102 & 19 &  & $<10^{-3}$ &  & 0.21 & 0.15 & 0.06 &  & 3.0 & 4.8 & 7.8 &  \\
\cmidrule(lr){4-10}\cmidrule(lr){11-25}
 &  & \multirow[c]{3}{*}{2} & \multirow[c]{3}{*}{776} & \multirow[c]{3}{*}{3320} & \multirow[c]{3}{*}{16} & \multirow[c]{3}{*}{393} & \multirow[c]{3}{*}{1169} & \multirow[c]{3}{*}{52.32} & \multirow[c]{3}{*}{6.27} &  & 10 & 410 & 78 &  & $<10^{-3}$ &  & 0.73 & 0.02 & 0.71 &  & 1.4 & 11.2 & 12.6 &  \\
 &  &  &  &  &  &  &  &  &  &  & \textbf{25} & 164 & 31 &  & $<10^{-3}$ &  & 0.14 & 0.01 & 0.13 &  & 3.0 & 4.4 & 7.4 &  \\
 &  &  &  &  &  &  &  &  &  &  & 40 & 102 & 19 &  & $<10^{-3}$ &  & 0.14 & 0.11 & 0.03 &  & 3.0 & 4.3 & 7.3 &  \\
\cmidrule(lr){4-10}\cmidrule(lr){11-25}
 &  & \multirow[c]{3}{*}{3} & \multirow[c]{3}{*}{776} & \multirow[c]{3}{*}{3320} & \multirow[c]{3}{*}{16} & \multirow[c]{3}{*}{439} & \multirow[c]{3}{*}{1215} & \multirow[c]{3}{*}{55.03} & \multirow[c]{3}{*}{6.55} &  & 10 & 410 & 78 &  & $<10^{-3}$ &  & 0.68 & 0.00 & 0.68 &  & 1.4 & 8.4 & 9.8 &  \\
 &  &  &  &  &  &  &  &  &  &  & 25 & 164 & 31 &  & $<10^{-3}$ &  & 0.21 & 0.06 & 0.15 &  & 3.0 & 3.9 & 6.9 &  \\
 &  &  &  &  &  &  &  &  &  &  & \textbf{40} & 102 & 19 &  & $<10^{-3}$ &  & 0.06 & 0.03 & 0.03 &  & 3.0 & 5.4 & 8.4 &  \\
\addlinespace[1pt]
\midrule
\multirow[c]{18}{*}{\shortstack{\texttt{Mallat}\\\href{https://www.timeseriesclassification.com/}{UCR}}} & \multirow[c]{9}{*}{20} & \multirow[c]{3}{*}{1} & \multirow[c]{3}{*}{360} & \multirow[c]{3}{*}{2040} & \multirow[c]{3}{*}{9} & \multirow[c]{3}{*}{240} & \multirow[c]{3}{*}{600} & \multirow[c]{3}{*}{25.56} & \multirow[c]{3}{*}{4.29} &  & 6 & 400 & 60 &  & 18.074 &  & 0.72 & 0.03 & 0.69 &  & 3.0 & 3.6 & 6.6 &  \\
 &  &  &  &  &  &  &  &  &  &  & \textbf{15} & 160 & 24 &  & 111.395 &  & 0.26 & 0.16 & 0.10 &  & 3.0 & 4.1 & 7.1 &  \\
 &  &  &  &  &  &  &  &  &  &  & 24 & 100 & 15 &  & 265.579 &  & 0.51 & 0.41 & 0.10 &  & 2.3 & 3.1 & 5.4 &  \\
\cmidrule(lr){4-10}\cmidrule(lr){11-25}
 &  & \multirow[c]{3}{*}{2} & \multirow[c]{3}{*}{360} & \multirow[c]{3}{*}{2040} & \multirow[c]{3}{*}{9} & \multirow[c]{3}{*}{245} & \multirow[c]{3}{*}{605} & \multirow[c]{3}{*}{26.39} & \multirow[c]{3}{*}{4.55} &  & 6 & 400 & 60 &  & 18.979 &  & 0.74 & 0.02 & 0.72 &  & 3.0 & 3.6 & 6.6 &  \\
 &  &  &  &  &  &  &  &  &  &  & \textbf{15} & 160 & 24 &  & 112.623 &  & 0.26 & 0.15 & 0.11 &  & 3.0 & 2.3 & 5.3 &  \\
 &  &  &  &  &  &  &  &  &  &  & 24 & 100 & 15 &  & 283.677 &  & 0.44 & 0.39 & 0.05 &  & 1.5 & 2.2 & 3.7 &  \\
\cmidrule(lr){4-10}\cmidrule(lr){11-25}
 &  & \multirow[c]{3}{*}{3} & \multirow[c]{3}{*}{360} & \multirow[c]{3}{*}{2040} & \multirow[c]{3}{*}{8} & \multirow[c]{3}{*}{226} & \multirow[c]{3}{*}{586} & \multirow[c]{3}{*}{25.00} & \multirow[c]{3}{*}{4.35} &  & 6 & 400 & 60 &  & 17.734 &  & 0.72 & 0.03 & 0.69 &  & 3.0 & 6.0 & 9.0 &  \\
 &  &  &  &  &  &  &  &  &  &  & \textbf{15} & 160 & 24 &  & 110.946 &  & 0.30 & 0.16 & 0.14 &  & 3.0 & 3.1 & 6.1 &  \\
 &  &  &  &  &  &  &  &  &  &  & 24 & 100 & 15 &  & 299.046 &  & 0.57 & 0.43 & 0.14 &  & 2.3 & 38.4 & 40.7 &  \\
\addlinespace[1pt]
\cmidrule(lr){3-10}\cmidrule(lr){11-25}
 & \multirow[c]{9}{*}{40} & \multirow[c]{3}{*}{1} & \multirow[c]{3}{*}{726} & \multirow[c]{3}{*}{1674} & \multirow[c]{3}{*}{72} & \multirow[c]{3}{*}{253} & \multirow[c]{3}{*}{979} & \multirow[c]{3}{*}{23.97} & \multirow[c]{3}{*}{7.63} &  & 6 & 400 & 121 &  & 25.063 &  & 1.54 & 0.03 & 1.51 &  & 3.0 & 2.7 & 5.7 &  \\
 &  &  &  &  &  &  &  &  &  &  & \textbf{13} & 185 & 56 &  & 116.515 &  & 0.54 & 0.20 & 0.34 &  & 3.0 & 2.4 & 5.4 &  \\
 &  &  &  &  &  &  &  &  &  &  & 21 & 114 & 35 &  & 428.516 &  & 0.76 & 0.62 & 0.14 &  & 2.3 & 2.9 & 5.2 &  \\
\cmidrule(lr){4-10}\cmidrule(lr){11-25}
 &  & \multirow[c]{3}{*}{2} & \multirow[c]{3}{*}{726} & \multirow[c]{3}{*}{1674} & \multirow[c]{3}{*}{104} & \multirow[c]{3}{*}{289} & \multirow[c]{3}{*}{1015} & \multirow[c]{3}{*}{23.28} & \multirow[c]{3}{*}{7.60} &  & 6 & 400 & 121 &  & 23.851 &  & 1.27 & 0.03 & 1.24 &  & 3.0 & 3.1 & 6.1 &  \\
 &  &  &  &  &  &  &  &  &  &  & \textbf{12} & 200 & 61 &  & 95.104 &  & 0.56 & 0.13 & 0.43 &  & 3.0 & 3.5 & 6.5 &  \\
 &  &  &  &  &  &  &  &  &  &  & 18 & 133 & 40 &  & 275.827 &  & 0.57 & 0.42 & 0.15 &  & 2.2 & 2.5 & 4.7 &  \\
\cmidrule(lr){4-10}\cmidrule(lr){11-25}
 &  & \multirow[c]{3}{*}{3} & \multirow[c]{3}{*}{726} & \multirow[c]{3}{*}{1674} & \multirow[c]{3}{*}{80} & \multirow[c]{3}{*}{259} & \multirow[c]{3}{*}{985} & \multirow[c]{3}{*}{25.21} & \multirow[c]{3}{*}{8.13} &  & 6 & 400 & 121 &  & 24.588 &  & 1.46 & 0.03 & 1.43 &  & 3.0 & 3.1 & 6.1 &  \\
 &  &  &  &  &  &  &  &  &  &  & \textbf{12} & 200 & 61 &  & 99.579 &  & 0.51 & 0.15 & 0.36 &  & 3.0 & 2.5 & 5.6 &  \\
 &  &  &  &  &  &  &  &  &  &  & 19 & 126 & 38 &  & 258.607 &  & 0.59 & 0.39 & 0.20 &  & 3.1 & 2.8 & 5.9 &  \\
\addlinespace[1pt]
\pagebreak
\multirow[c]{18}{*}{\shortstack{\texttt{PR2392}\\\href{https://www.math.uwaterloo.ca/tsp/history/tspinfo/pr2392_info.html}{TSPLIB}}} & \multirow[c]{9}{*}{20} & \multirow[c]{3}{*}{1} & \multirow[c]{3}{*}{387} & \multirow[c]{3}{*}{2005} & \multirow[c]{3}{*}{12} & \multirow[c]{3}{*}{222} & \multirow[c]{3}{*}{609} & \multirow[c]{3}{*}{31.78} & \multirow[c]{3}{*}{3.53} &  & 5 & 478 & 77 &  & $<10^{-3}$ &  & 4.96 & 0.00 & 4.96 &  & 0.3 & 6.5 & 6.8 &  \\
 &  &  &  &  &  &  &  &  &  &  & \textbf{14} & 171 & 28 &  & 0.001 &  & 1.46 & 0.06 & 1.40 &  & 1.7 & 2.3 & 4.0 &  \\
 &  &  &  &  &  &  &  &  &  &  & 23 & 104 & 17 &  & 0.004 &  & 2.18 & 1.37 & 0.81 &  & 1.1 & 4.5 & 5.5 &  \\
\cmidrule(lr){4-10}\cmidrule(lr){11-25}
 &  & \multirow[c]{3}{*}{2} & \multirow[c]{3}{*}{387} & \multirow[c]{3}{*}{2005} & \multirow[c]{3}{*}{16} & \multirow[c]{3}{*}{209} & \multirow[c]{3}{*}{596} & \multirow[c]{3}{*}{30.75} & \multirow[c]{3}{*}{3.23} &  & 5 & 478 & 77 &  & $<10^{-3}$ &  & 4.10 & 0.00 & 4.10 &  & 0.3 & 7.2 & 7.5 &  \\
 &  &  &  &  &  &  &  &  &  &  & 14 & 171 & 28 &  & 0.001 &  & 1.45 & 0.13 & 1.32 &  & 1.7 & 2.4 & 4.0 &  \\
 &  &  &  &  &  &  &  &  &  &  & \textbf{23} & 104 & 17 &  & 0.004 &  & 1.28 & 0.42 & 0.86 &  & 1.1 & 4.3 & 5.4 &  \\
\cmidrule(lr){4-10}\cmidrule(lr){11-25}
 &  & \multirow[c]{3}{*}{3} & \multirow[c]{3}{*}{387} & \multirow[c]{3}{*}{2005} & \multirow[c]{3}{*}{24} & \multirow[c]{3}{*}{209} & \multirow[c]{3}{*}{596} & \multirow[c]{3}{*}{30.75} & \multirow[c]{3}{*}{3.68} &  & 5 & 478 & 77 &  & $<10^{-3}$ &  & 4.50 & 0.01 & 4.49 &  & 0.3 & 7.7 & 8.0 &  \\
 &  &  &  &  &  &  &  &  &  &  & 14 & 171 & 28 &  & 0.001 &  & 2.39 & 0.88 & 1.51 &  & 0.7 & 2.4 & 3.1 &  \\
 &  &  &  &  &  &  &  &  &  &  & \textbf{23} & 104 & 17 &  & 0.004 &  & 1.78 & 1.17 & 0.61 &  & 2.7 & 9.4 & 12.1 &  \\
\addlinespace[1pt]
\cmidrule(lr){3-10}\cmidrule(lr){11-25}
 & \multirow[c]{9}{*}{40} & \multirow[c]{3}{*}{1} & \multirow[c]{3}{*}{774} & \multirow[c]{3}{*}{1620} & \multirow[c]{3}{*}{103} & \multirow[c]{3}{*}{228} & \multirow[c]{3}{*}{1000} & \multirow[c]{3}{*}{31.78} & \multirow[c]{3}{*}{7.09} &  & 5 & 478 & 155 &  & $<10^{-3}$ &  & 3.99 & 0.02 & 3.97 &  & 0.2 & 5.8 & 6.1 &  \\
 &  &  &  &  &  &  &  &  &  &  & \textbf{13} & 184 & 60 &  & 0.001 &  & 1.31 & 0.10 & 1.21 &  & 1.3 & 2.8 & 4.0 &  \\
 &  &  &  &  &  &  &  &  &  &  & 21 & 114 & 37 &  & 37.960 &  & 1.75 & 0.71 & 1.04 &  & 1.9 & 8.8 & 10.7 &  \\
\cmidrule(lr){4-10}\cmidrule(lr){11-25}
 &  & \multirow[c]{3}{*}{2} & \multirow[c]{3}{*}{774} & \multirow[c]{3}{*}{1619} & \multirow[c]{3}{*}{223} & \multirow[c]{3}{*}{228} & \multirow[c]{3}{*}{1001} & \multirow[c]{3}{*}{34.37} & \multirow[c]{3}{*}{6.94} &  & 5 & 478 & 155 &  & $<10^{-3}$ &  & 3.96 & 0.01 & 3.95 &  & 0.2 & 4.8 & 5.0 &  \\
 &  &  &  &  &  &  &  &  &  &  & 7 & 342 & 111 &  & 0.133 &  & 2.62 & 0.01 & 2.61 &  & 0.3 & 4.8 & 5.1 &  \\
 &  &  &  &  &  &  &  &  &  &  & \textbf{10} & 239 & 77 &  & 131.818 &  & 2.14 & 0.72 & 1.42 &  & 0.4 & 12.9 & 13.3 &  \\
\cmidrule(lr){4-10}\cmidrule(lr){11-25}
 &  & \multirow[c]{3}{*}{3} & \multirow[c]{3}{*}{774} & \multirow[c]{3}{*}{1618} & \multirow[c]{3}{*}{140} & \multirow[c]{3}{*}{228} & \multirow[c]{3}{*}{1002} & \multirow[c]{3}{*}{36.30} & \multirow[c]{3}{*}{7.15} &  & 5 & 478 & 155 &  & $<10^{-3}$ &  & 3.91 & 0.02 & 3.89 &  & 0.2 & 6.3 & 6.6 &  \\
 &  &  &  &  &  &  &  &  &  &  & \textbf{10} & 239 & 77 &  & 0.001 &  & 1.88 & 0.01 & 1.87 &  & 1.0 & 3.5 & 4.5 &  \\
 &  &  &  &  &  &  &  &  &  &  & 16 & 150 & 48 &  & 89.693 &  & 1.91 & 0.59 & 1.32 &  & 2.3 & 4.5 & 6.9 &  \\
\addlinespace[1pt]
\midrule
\multirow[c]{18}{*}{\shortstack{\texttt{Rice}\\\href{https://archive.ics.uci.edu/dataset/545/rice}{UCI}}} & \multirow[c]{9}{*}{20} & \multirow[c]{3}{*}{1} & \multirow[c]{3}{*}{315} & \multirow[c]{3}{*}{3495} & \multirow[c]{3}{*}{7} & \multirow[c]{3}{*}{244} & \multirow[c]{3}{*}{559} & \multirow[c]{3}{*}{58.10} & \multirow[c]{3}{*}{3.51} &  & 9 & 423 & 35 &  & 0.007 &  & 0.24 & 0.00 & 0.24 &  & 2.7 & 8.1 & 10.8 &  \\
 &  &  &  &  &  &  &  &  &  &  & 23 & 166 & 14 &  & 0.087 &  & 0.07 & 0.00 & 0.07 &  & 3.0 & 5.1 & 8.1 &  \\
 &  &  &  &  &  &  &  &  &  &  & \textbf{38} & 100 & 8 &  & 0.204 &  & 0.04 & 0.02 & 0.02 &  & 3.0 & 4.1 & 7.1 &  \\
\cmidrule(lr){4-10}\cmidrule(lr){11-25}
 &  & \multirow[c]{3}{*}{2} & \multirow[c]{3}{*}{315} & \multirow[c]{3}{*}{3495} & \multirow[c]{3}{*}{6} & \multirow[c]{3}{*}{242} & \multirow[c]{3}{*}{557} & \multirow[c]{3}{*}{57.78} & \multirow[c]{3}{*}{3.40} &  & 9 & 423 & 35 &  & 0.015 &  & 0.38 & 0.00 & 0.38 &  & 2.7 & 7.2 & 9.9 &  \\
 &  &  &  &  &  &  &  &  &  &  & 23 & 166 & 14 &  & 0.094 &  & 0.07 & 0.01 & 0.06 &  & 3.0 & 4.6 & 7.6 &  \\
 &  &  &  &  &  &  &  &  &  &  & \textbf{38} & 100 & 8 &  & 0.184 &  & 0.04 & 0.03 & 0.01 &  & 3.0 & 4.1 & 7.1 &  \\
\cmidrule(lr){4-10}\cmidrule(lr){11-25}
 &  & \multirow[c]{3}{*}{3} & \multirow[c]{3}{*}{315} & \multirow[c]{3}{*}{3495} & \multirow[c]{3}{*}{7} & \multirow[c]{3}{*}{232} & \multirow[c]{3}{*}{547} & \multirow[c]{3}{*}{55.24} & \multirow[c]{3}{*}{3.30} &  & 9 & 423 & 35 &  & 0.020 &  & 0.37 & 0.00 & 0.37 &  & 2.7 & 6.8 & 9.5 &  \\
 &  &  &  &  &  &  &  &  &  &  & 23 & 166 & 14 &  & 0.070 &  & 0.05 & 0.01 & 0.04 &  & 3.0 & 3.6 & 6.6 &  \\
 &  &  &  &  &  &  &  &  &  &  & \textbf{38} & 100 & 8 &  & 0.189 &  & 0.04 & 0.03 & 0.01 &  & 3.0 & 3.6 & 6.6 &  \\
\addlinespace[1pt]
\cmidrule(lr){3-10}\cmidrule(lr){11-25}
 & \multirow[c]{9}{*}{40} & \multirow[c]{3}{*}{1} & \multirow[c]{3}{*}{630} & \multirow[c]{3}{*}{3180} & \multirow[c]{3}{*}{16} & \multirow[c]{3}{*}{357} & \multirow[c]{3}{*}{987} & \multirow[c]{3}{*}{57.30} & \multirow[c]{3}{*}{6.22} &  & 9 & 423 & 70 &  & 0.018 &  & 0.15 & 0.00 & 0.15 &  & 2.6 & 6.9 & 9.5 &  \\
 &  &  &  &  &  &  &  &  &  &  & 23 & 166 & 27 &  & 0.105 &  & 0.03 & 0.00 & 0.03 &  & 3.0 & 33.8 & 36.8 &  \\
 &  &  &  &  &  &  &  &  &  &  & \textbf{38} & 100 & 17 &  & 0.197 &  & 0.01 & 0.01 & 0.00 &  & 3.0 & 2.6 & 5.7 &  \\
\cmidrule(lr){4-10}\cmidrule(lr){11-25}
 &  & \multirow[c]{3}{*}{2} & \multirow[c]{3}{*}{630} & \multirow[c]{3}{*}{3180} & \multirow[c]{3}{*}{10} & \multirow[c]{3}{*}{357} & \multirow[c]{3}{*}{987} & \multirow[c]{3}{*}{55.08} & \multirow[c]{3}{*}{6.70} &  & 9 & 423 & 70 &  & 0.016 &  & 0.25 & 0.00 & 0.25 &  & 2.2 & 6.3 & 8.5 &  \\
 &  &  &  &  &  &  &  &  &  &  & \textbf{23} & 166 & 27 &  & 0.101 &  & 0.02 & 0.00 & 0.02 &  & 3.1 & 3.8 & 6.9 &  \\
 &  &  &  &  &  &  &  &  &  &  & 38 & 100 & 17 &  & 0.227 &  & 0.02 & 0.02 & 0.00 &  & 3.0 & 2.9 & 5.9 &  \\
\cmidrule(lr){4-10}\cmidrule(lr){11-25}
 &  & \multirow[c]{3}{*}{3} & \multirow[c]{3}{*}{630} & \multirow[c]{3}{*}{3180} & \multirow[c]{3}{*}{13} & \multirow[c]{3}{*}{368} & \multirow[c]{3}{*}{998} & \multirow[c]{3}{*}{59.37} & \multirow[c]{3}{*}{6.56} &  & 9 & 423 & 70 &  & 0.015 &  & 0.33 & 0.00 & 0.33 &  & 2.3 & 5.9 & 8.3 &  \\
 &  &  &  &  &  &  &  &  &  &  & 23 & 166 & 27 &  & 0.134 &  & 0.05 & 0.00 & 0.05 &  & 1.8 & 45.0 & 46.8 &  \\
 &  &  &  &  &  &  &  &  &  &  & \textbf{38} & 100 & 17 &  & 0.219 &  & 0.01 & 0.01 & 0.00 &  & 3.0 & 3.1 & 6.1 &  \\
\addlinespace[1pt]
\pagebreak
\multirow[c]{18}{*}{\shortstack{\texttt{Wifi}\\\href{https://archive.ics.uci.edu/dataset/422/wireless+indoor+localization}{UCI}}} & \multirow[c]{9}{*}{20} & \multirow[c]{3}{*}{1} & \multirow[c]{3}{*}{130} & \multirow[c]{3}{*}{1870} & \multirow[c]{3}{*}{5} & \multirow[c]{3}{*}{107} & \multirow[c]{3}{*}{237} & \multirow[c]{3}{*}{24.62} & \multirow[c]{3}{*}{1.96} &  & 5 & 400 & 26 &  & 9.020 &  & 0.33 & 0.01 & 0.32 &  & 1.0 & 4.0 & 5.0 &  \\
 &  &  &  &  &  &  &  &  &  &  & \textbf{11} & 182 & 12 &  & 64.324 &  & 0.19 & 0.13 & 0.06 &  & 2.9 & 17.9 & 20.8 &  \\
 &  &  &  &  &  &  &  &  &  &  & 17 & 118 & 8 &  & 325.048 &  & 0.42 & 0.39 & 0.03 &  & 1.6 & 1.8 & 3.4 &  \\
\cmidrule(lr){4-10}\cmidrule(lr){11-25}
 &  & \multirow[c]{3}{*}{2} & \multirow[c]{3}{*}{130} & \multirow[c]{3}{*}{1870} & \multirow[c]{3}{*}{5} & \multirow[c]{3}{*}{106} & \multirow[c]{3}{*}{236} & \multirow[c]{3}{*}{20.77} & \multirow[c]{3}{*}{1.68} &  & 5 & 400 & 26 &  & 9.840 &  & 0.26 & 0.04 & 0.22 &  & 1.0 & 4.0 & 5.1 &  \\
 &  &  &  &  &  &  &  &  &  &  & \textbf{11} & 182 & 12 &  & 61.682 &  & 0.15 & 0.10 & 0.05 &  & 1.3 & 2.2 & 3.6 &  \\
 &  &  &  &  &  &  &  &  &  &  & 17 & 118 & 8 &  & 398.545 &  & 0.43 & 0.39 & 0.04 &  & 1.6 & 1.6 & 3.2 &  \\
\cmidrule(lr){4-10}\cmidrule(lr){11-25}
 &  & \multirow[c]{3}{*}{3} & \multirow[c]{3}{*}{130} & \multirow[c]{3}{*}{1870} & \multirow[c]{3}{*}{5} & \multirow[c]{3}{*}{96} & \multirow[c]{3}{*}{226} & \multirow[c]{3}{*}{22.31} & \multirow[c]{3}{*}{1.87} &  & 5 & 400 & 26 &  & 8.720 &  & 0.23 & 0.01 & 0.22 &  & 1.0 & 5.3 & 6.3 &  \\
 &  &  &  &  &  &  &  &  &  &  & \textbf{9} & 222 & 14 &  & 47.982 &  & 0.17 & 0.07 & 0.10 &  & 1.2 & 2.2 & 3.4 &  \\
 &  &  &  &  &  &  &  &  &  &  & 13 & 154 & 10 &  & 137.541 &  & 0.23 & 0.20 & 0.03 &  & 1.4 & 1.7 & 3.1 &  \\
\addlinespace[1pt]
\cmidrule(lr){3-10}\cmidrule(lr){11-25}
 & \multirow[c]{9}{*}{40} & \multirow[c]{3}{*}{1} & \multirow[c]{3}{*}{260} & \multirow[c]{3}{*}{1740} & \multirow[c]{3}{*}{8} & \multirow[c]{3}{*}{160} & \multirow[c]{3}{*}{420} & \multirow[c]{3}{*}{21.15} & \multirow[c]{3}{*}{2.53} &  & 5 & 400 & 52 &  & 9.930 &  & 0.34 & 0.02 & 0.32 &  & 0.9 & 5.7 & 6.6 &  \\
 &  &  &  &  &  &  &  &  &  &  & \textbf{11} & 182 & 24 &  & 82.851 &  & 0.20 & 0.10 & 0.10 &  & 1.1 & 2.4 & 3.5 &  \\
 &  &  &  &  &  &  &  &  &  &  & 17 & 118 & 15 &  & 335.933 &  & 0.32 & 0.29 & 0.03 &  & 1.4 & 1.5 & 3.0 &  \\
\cmidrule(lr){4-10}\cmidrule(lr){11-25}
 &  & \multirow[c]{3}{*}{2} & \multirow[c]{3}{*}{260} & \multirow[c]{3}{*}{1740} & \multirow[c]{3}{*}{9} & \multirow[c]{3}{*}{170} & \multirow[c]{3}{*}{430} & \multirow[c]{3}{*}{21.92} & \multirow[c]{3}{*}{3.09} &  & 5 & 400 & 52 &  & 6.710 &  & 0.39 & 0.02 & 0.37 &  & 0.9 & 6.2 & 7.1 &  \\
 &  &  &  &  &  &  &  &  &  &  & 9 & 222 & 29 &  & 37.402 &  & 0.25 & 0.07 & 0.18 &  & 1.1 & 2.9 & 4.0 &  \\
 &  &  &  &  &  &  &  &  &  &  & \textbf{13} & 154 & 20 &  & 126.728 &  & 0.23 & 0.16 & 0.07 &  & 1.3 & 2.1 & 3.4 &  \\
\cmidrule(lr){4-10}\cmidrule(lr){11-25}
 &  & \multirow[c]{3}{*}{3} & \multirow[c]{3}{*}{260} & \multirow[c]{3}{*}{1740} & \multirow[c]{3}{*}{8} & \multirow[c]{3}{*}{176} & \multirow[c]{3}{*}{436} & \multirow[c]{3}{*}{26.54} & \multirow[c]{3}{*}{3.63} &  & 5 & 400 & 52 &  & 9.970 &  & 0.43 & 0.02 & 0.41 &  & 1.0 & 3.9 & 4.9 &  \\
 &  &  &  &  &  &  &  &  &  &  & \textbf{9} & 222 & 29 &  & 45.100 &  & 0.14 & 0.07 & 0.07 &  & 1.2 & 3.5 & 4.7 &  \\
 &  &  &  &  &  &  &  &  &  &  & 13 & 154 & 20 &  & 117.549 &  & 0.25 & 0.15 & 0.10 &  & 1.3 & 1.6 & 3.0 &  \\
\addlinespace[1pt]
\end{longtable}
\endgroup

\newpage
\setcounter{table}{0}
\section{Comparison between CAP, Gurobi and CN}
\label{appendixB}

Table~\ref{tab:partition_treatment_results} reports the results obtained with
CAP, the MISOCP solve using Gurobi (GRB), and the constraint-neglecting
variant CN for every instance. For each run, the table gives the selected
number of anticlusters $T^*$, the percentage of CL constraints retained by
CN, and the values of $D(\mathcal P,\mathcal S)$ obtained by the three
methods.
For GRB, a superscript ``$*$'' indicates an out-of-memory termination
for which Gurobi returns an incumbent, while ``$-$'' indicates that no
incumbent is found. In the GRB time column, ``$t\ell$'' denotes
termination at the time limit. The partitioning-time columns report the
running times of VNS, the GRB solve, and CN. The remaining columns
report the lower-bound computation times and the resulting certified gaps
for CAP and CN.

\begingroup
\small
\setlength{\tabcolsep}{1.2pt}
\renewcommand{\arraystretch}{0.9}
\begin{longtable}{lrrrr>{\hspace{2pt}}c<{\hspace{2pt}}rrr>{\hspace{2pt}}c<{\hspace{2pt}}rrr>{\hspace{2pt}}c<{\hspace{2pt}}rr>{\hspace{2pt}}c<{\hspace{2pt}}rr}
\caption{CAP, GRB, and CN detailed results.}
\label{tab:partition_treatment_results} \\
\toprule
Dataset & $p$ & Rep. & $T$ & $CL$ (\%) & \multicolumn{1}{c}{\hspace{6pt}} & \multicolumn{3}{c}{$D(\mathcal{P},\mathcal{S})$} & \multicolumn{1}{c}{\hspace{6pt}} & \multicolumn{3}{c}{Time $D(\mathcal P, \mathcal S)$} & \multicolumn{1}{c}{\hspace{6pt}} & \multicolumn{2}{c}{Time LB (min)} & \multicolumn{1}{c}{\hspace{6pt}} & \multicolumn{2}{c}{$\widehat\gamma$} \\
\cmidrule(lr){7-9} \cmidrule(lr){11-13} \cmidrule(lr){15-16} \cmidrule(lr){18-19}
 & & & & & \multicolumn{1}{c}{\hspace{6pt}} & CAP & GRB & CN & \multicolumn{1}{c}{\hspace{6pt}} & VNS & GRB & CN & \multicolumn{1}{c}{\hspace{6pt}} & CAP & CN & \multicolumn{1}{c}{\hspace{6pt}} & CAP & CN \\
\midrule
\endfirsthead
\multicolumn{19}{l}{\tablename~\thetable{} (continued)} \\
\toprule
Dataset & $p$ & Rep. & $T$ & $CL$ (\%) & \multicolumn{1}{c}{\hspace{6pt}} & \multicolumn{3}{c}{$D(\mathcal{P},\mathcal{S})$} & \multicolumn{1}{c}{\hspace{6pt}} & \multicolumn{3}{c}{Time $D(\mathcal P, \mathcal S)$} & \multicolumn{1}{c}{\hspace{6pt}} & \multicolumn{2}{c}{Time LB (min)} & \multicolumn{1}{c}{\hspace{6pt}} & \multicolumn{2}{c}{$\widehat\gamma$} \\
\cmidrule(lr){7-9} \cmidrule(lr){11-13} \cmidrule(lr){15-16} \cmidrule(lr){18-19}
 & & & & & \multicolumn{1}{c}{\hspace{6pt}} & CAP & GRB & CN & \multicolumn{1}{c}{\hspace{6pt}} & VNS & GRB & CN & \multicolumn{1}{c}{\hspace{6pt}} & CAP & CN & \multicolumn{1}{c}{\hspace{6pt}} & CAP & CN \\
\midrule
\endhead
\midrule
\multicolumn{19}{r}{continued on next page} \\
\endfoot
\bottomrule
\endlastfoot
\multirow[c]{6}{*}{\texttt{syn-1500-50}} & \multirow[c]{3}{*}{20} & 1 & 15 & 5.67 &  & $<10^{-3}$ & 0.987 & $<10^{-3}$ &  & 1.0 & $t\ell$ & 1.4 &  & 1.4 & 0.6 &  & 0.10 & 37.07 \\
 &  & 2 & 15 & 4.00 &  & $<10^{-3}$ & 1.938 & $<10^{-3}$ &  & 0.8 & $t\ell$ & 1.3 &  & 1.1 & 0.6 &  & 0.49 & 31.92 \\
 &  & 3 & 15 & 5.67 &  & $<10^{-3}$ & 0.479 & $<10^{-3}$ &  & 1.0 & $t\ell$ & 1.4 &  & 1.1 & 0.4 &  & 0.27 & 38.00 \\
\cmidrule(lr){3-19}
 & \multirow[c]{3}{*}{40} & 1 & 8 & 13.17 &  & 7.901 & 8.412 & $<10^{-3}$ &  & 0.3 & $t\ell$ & 0.4 &  & 1.1 & 0.8 &  & 3.94 & 45.13 \\
 &  & 2 & 11 & 9.50 &  & 7.873 & 74.557 & $<10^{-3}$ &  & 0.4 & $t\ell$ & 1.2 &  & 1.2 & 0.5 &  & 1.93 & 45.49 \\
 &  & 3 & 15 & 6.50 &  & 0.474 & 22.899 & $<10^{-3}$ &  & 0.5 & $t\ell$ & 1.4 &  & 1.5 & 0.6 &  & 0.40 & 46.93 \\
\midrule
\multirow[c]{6}{*}{\texttt{syn-1500-75}} & \multirow[c]{3}{*}{20} & 1 & 15 & 5.72 &  & $<10^{-3}$ & 2.268 & $<10^{-3}$ &  & 0.9 & $t\ell$ & 16.7 &  & 1.4 & 1.4 &  & 0.22 & 14.25 \\
 &  & 2 & 15 & 7.07 &  & $<10^{-3}$ & 1.435 & $<10^{-3}$ &  & 0.9 & $t\ell$ & 1.4 &  & 1.3 & 1.0 &  & 0.19 & 14.24 \\
 &  & 3 & 15 & 7.74 &  & $<10^{-3}$ & 1.131 & $<10^{-3}$ &  & 1.0 & $t\ell$ & 1.3 &  & 1.3 & 1.1 &  & 0.12 & 13.56 \\
\cmidrule(lr){3-19}
 & \multirow[c]{3}{*}{40} & 1 & 10 & 9.09 &  & 6.952 & 20.525 & $<10^{-3}$ &  & 0.5 & $t\ell$ & 1.0 &  & 1.0 & 0.9 &  & 1.59 & 22.41 \\
 &  & 2 & 15 & 4.88 &  & 0.001 & 23.622 & $<10^{-3}$ &  & 0.5 & $t\ell$ & 1.4 &  & 2.1 & 1.2 &  & 0.53 & 22.71 \\
 &  & 3 & 8 & 9.60 &  & 2.749 & 4.048 & $<10^{-3}$ &  & 0.4 & $t\ell$ & 1.0 &  & 1.2 & 0.7 &  & 1.66 & 21.50 \\
\midrule
\multirow[c]{6}{*}{\texttt{syn-3000-50}} & \multirow[c]{3}{*}{20} & 1 & 30 & 2.85 &  & $<10^{-3}$ & $-$ & $<10^{-3}$ &  & 3.2 & $t\ell$ & 3.2 &  & 2.4 & 0.9 &  & 0.05 & 34.31 \\
 &  & 2 & 30 & 2.85 &  & $<10^{-3}$ & $-$ & $<10^{-3}$ &  & 3.1 & $t\ell$ & 3.1 &  & 2.4 & 0.9 &  & 0.17 & 34.83 \\
 &  & 3 & 30 & 3.18 &  & $<10^{-3}$ & $-$ & $<10^{-3}$ &  & 3.1 & $t\ell$ & 3.0 &  & 2.7 & 1.1 &  & 0.29 & 32.65 \\
\cmidrule(lr){3-19}
 & \multirow[c]{3}{*}{40} & 1 & 21 & 5.43 &  & 3.957 & $-$ & $<10^{-3}$ &  & 3.0 & $t\ell$ & 3.0 &  & 2.0 & 0.9 &  & 1.06 & 44.96 \\
 &  & 2 & 18 & 5.93 &  & 22.281 & 31.253 & $<10^{-3}$ &  & 2.5 & $t\ell$ & 3.0 &  & 2.1 & 0.9 &  & 3.24 & 46.41 \\
 &  & 3 & 12 & 8.44 &  & 18.300 & 25.538 & $<10^{-3}$ &  & 2.0 & $t\ell$ & 3.0 &  & 3.0 & 5.5 &  & 4.18 & 44.01 \\
\midrule
\multirow[c]{6}{*}{\texttt{syn-3000-75}} & \multirow[c]{3}{*}{20} & 1 & 30 & 4.06 &  & 0.001 & 230.751 & $<10^{-3}$ &  & 2.3 & $t\ell$ & 3.1 &  & 2.5 & 5.4 &  & 0.20 & 12.96 \\
 &  & 2 & 30 & 3.38 &  & 0.001 & 238.594 & $<10^{-3}$ &  & 3.0 & $t\ell$ & 3.1 &  & 2.8 & 13.8 &  & 0.23 & 13.11 \\
 &  & 3 & 30 & 3.38 &  & $<10^{-3}$ & $-$ & $<10^{-3}$ &  & 3.0 & $t\ell$ & 3.0 &  & 3.1 & 1.7 &  & 0.37 & 15.05 \\
\cmidrule(lr){3-19}
 & \multirow[c]{3}{*}{40} & 1 & 6 & 16.67 &  & 3.700 & 19.516 & $<10^{-3}$ &  & 0.6 & $t\ell$ & 3.0 &  & 4.4 & 2.4 &  & 3.53 & 18.67 \\
 &  & 2 & 4 & 24.28 &  & 0.010 & 0.307 & $<10^{-3}$ &  & 0.4 & $t\ell$ & 2.7 &  & 8.9 & 5.5 &  & 4.79 & 16.77 \\
 &  & 3 & 6 & 17.85 &  & 11.432 & 66.086 & $<10^{-3}$ &  & 1.2 & $t\ell$ & 3.0 &  & 5.4 & 7.6 &  & 3.52 & 18.48 \\
\pagebreak
\multirow[c]{6}{*}{\texttt{Abalone}} & \multirow[c]{3}{*}{20} & 1 & 39 & 2.94 &  & 1.101 & 17.785\textsuperscript{*} & 1.049 &  & 2.8 & 9.1 & 3.8 &  & 4.1 & 1.8 &  & 0.76 & 23.02 \\
 &  & 2 & 40 & 2.81 &  & 1.101 & 19.624\textsuperscript{*} & 1.096 &  & 2.8 & 8.8 & 2.3 &  & 4.8 & 2.2 &  & 0.74 & 23.15 \\
 &  & 3 & 39 & 2.04 &  & 1.115 & 15.271\textsuperscript{*} & 0.856 &  & 3.4 & 9.2 & 3.2 &  & 3.6 & 1.5 &  & 0.66 & 20.07 \\
\cmidrule(lr){3-19}
 & \multirow[c]{3}{*}{40} & 1 & 31 & 3.12 &  & 5.804 & $-$ & 0.442 &  & 1.5 & $t\ell$ & 3.0 &  & 2.6 & 1.6 &  & 1.93 & 31.06 \\
 &  & 2 & 28 & 3.82 &  & 2.352 & 9.738 & 0.358 &  & 3.4 & $t\ell$ & 3.0 &  & 2.6 & 1.6 &  & 1.81 & 30.51 \\
 &  & 3 & 34 & 2.80 &  & 3.722 & $-$ & 0.787 &  & 2.5 & 7.2 & 3.0 &  & 3.3 & 1.6 &  & 1.37 & 32.46 \\
\midrule
\multirow[c]{6}{*}{\texttt{EngyTime}} & \multirow[c]{3}{*}{20} & 1 & 40 & 2.32 &  & $<10^{-3}$ & $-$ & $<10^{-3}$ &  & 3.0 & 2.3 & 3.0 &  & 7.3 & 6.6 &  & 0.36 & 4.79 \\
 &  & 2 & 40 & 3.61 &  & $<10^{-3}$ & $-$ & $<10^{-3}$ &  & 3.0 & 2.1 & 3.0 &  & 6.4 & 7.7 &  & 0.29 & 4.41 \\
 &  & 3 & 25 & 4.38 &  & $<10^{-3}$ & 200.687\textsuperscript{*} & $<10^{-3}$ &  & 3.0 & 4.3 & 3.0 &  & 5.1 & 5.6 &  & 0.39 & 4.37 \\
\cmidrule(lr){3-19}
 & \multirow[c]{3}{*}{40} & 1 & 40 & 3.22 &  & $<10^{-3}$ & $-$ & $<10^{-3}$ &  & 3.0 & 1.9 & 3.0 &  & 4.8 & 6.5 &  & 0.21 & 7.85 \\
 &  & 2 & 25 & 5.03 &  & $<10^{-3}$ & 420.404\textsuperscript{*} & $<10^{-3}$ &  & 3.0 & 4.6 & 3.0 &  & 4.4 & 5.3 &  & 0.14 & 7.08 \\
 &  & 3 & 40 & 2.19 &  & $<10^{-3}$ & $-$ & $<10^{-3}$ &  & 3.0 & 1.8 & 3.0 &  & 5.4 & 6.7 &  & 0.06 & 7.69 \\
\midrule
\multirow[c]{6}{*}{\texttt{Mallat}} & \multirow[c]{3}{*}{20} & 1 & 15 & 6.67 &  & 111.395 & $-$ & 95.367 &  & 3.0 & $-$ & 3.0 &  & 4.1 & 1.3 &  & 0.26 & 4.21 \\
 &  & 2 & 15 & 5.00 &  & 112.623 & $-$ & 96.345 &  & 3.0 & $-$ & 3.0 &  & 2.3 & 1.4 &  & 0.26 & 4.60 \\
 &  & 3 & 15 & 5.28 &  & 110.946 & 431.241 & 98.395 &  & 3.0 & $t\ell$ & 3.0 &  & 3.1 & 1.2 &  & 0.30 & 4.47 \\
\cmidrule(lr){3-19}
 & \multirow[c]{3}{*}{40} & 1 & 13 & 7.16 &  & 116.515 & 382.174 & 74.828 &  & 3.0 & $t\ell$ & 3.0 &  & 2.4 & 1.4 &  & 0.54 & 7.39 \\
 &  & 2 & 12 & 10.74 &  & 95.104 & 353.468 & 64.713 &  & 3.0 & $t\ell$ & 3.0 &  & 3.5 & 2.0 &  & 0.56 & 7.08 \\
 &  & 3 & 12 & 8.82 &  & 99.579 & 323.072 & 58.260 &  & 3.0 & $t\ell$ & 3.0 &  & 2.5 & 1.9 &  & 0.51 & 7.49 \\
\midrule
\multirow[c]{6}{*}{\texttt{PR2392}} & \multirow[c]{3}{*}{20} & 1 & 14 & 6.46 &  & 0.001 & 0.266 & 0.001 &  & 1.7 & $t\ell$ & 0.8 &  & 2.3 & 2.6 &  & 1.46 & 7.24 \\
 &  & 2 & 23 & 2.33 &  & 0.004 & 18.085 & 0.004 &  & 1.1 & $t\ell$ & 1.3 &  & 4.3 & 4.4 &  & 1.28 & 6.94 \\
 &  & 3 & 23 & 6.20 &  & 0.004 & 32.374 & 0.003 &  & 2.7 & $t\ell$ & 1.2 &  & 9.4 & 3.8 &  & 1.78 & 6.91 \\
\cmidrule(lr){3-19}
 & \multirow[c]{3}{*}{40} & 1 & 13 & 8.40 &  & 0.001 & 0.291 & 0.001 &  & 1.3 & $t\ell$ & 0.8 &  & 2.8 & 2.5 &  & 1.31 & 9.15 \\
 &  & 2 & 10 & 10.21 &  & 131.818 & 136.329 & 0.001 &  & 0.4 & $t\ell$ & 0.7 &  & 12.9 & 4.0 &  & 2.14 & 9.69 \\
 &  & 3 & 10 & 9.04 &  & 0.001 & 0.181 & $<10^{-3}$ &  & 1.0 & $t\ell$ & 0.7 &  & 3.5 & 4.9 &  & 1.88 & 10.01 \\
\midrule
\multirow[c]{6}{*}{\texttt{Rice}} & \multirow[c]{3}{*}{20} & 1 & 38 & 3.81 &  & 0.204 & $-$ & 0.186 &  & 3.0 & 4.3 & 3.0 &  & 4.1 & 4.2 &  & 0.04 & 3.67 \\
 &  & 2 & 38 & 1.91 &  & 0.184 & $-$ & 0.219 &  & 3.0 & 4.5 & 3.0 &  & 4.1 & 5.0 &  & 0.04 & 3.75 \\
 &  & 3 & 38 & 0.95 &  & 0.189 & $-$ & 0.231 &  & 3.0 & 4.3 & 3.0 &  & 3.6 & 4.5 &  & 0.04 & 3.67 \\
\cmidrule(lr){3-19}
 & \multirow[c]{3}{*}{40} & 1 & 38 & 3.33 &  & 0.197 & $-$ & 0.194 &  & 3.0 & 4.4 & 3.0 &  & 2.6 & 4.3 &  & 0.01 & 6.64 \\
 &  & 2 & 23 & 3.18 &  & 0.101 & 847.098 & 0.109 &  & 3.1 & $t\ell$ & 3.0 &  & 3.8 & 3.6 &  & 0.02 & 6.87 \\
 &  & 3 & 38 & 2.54 &  & 0.219 & $-$ & 0.176 &  & 3.0 & 4.1 & 3.0 &  & 3.1 & 4.5 &  & 0.01 & 6.98 \\
\midrule
\multirow[c]{6}{*}{\texttt{Wifi}} & \multirow[c]{3}{*}{20} & 1 & 11 & 5.39 &  & 64.324 & 356.255 & 64.738 &  & 2.9 & $t\ell$ & 1.6 &  & 17.9 & 1.5 &  & 0.19 & 2.17 \\
 &  & 2 & 11 & 6.15 &  & 61.682 & 271.180 & 70.220 &  & 1.3 & $t\ell$ & 1.7 &  & 2.2 & 1.8 &  & 0.15 & 1.91 \\
 &  & 3 & 9 & 15.39 &  & 47.982 & 181.038 & 41.218 &  & 1.2 & $t\ell$ & 1.4 &  & 2.2 & 2.0 &  & 0.17 & 1.85 \\
\cmidrule(lr){3-19}
 & \multirow[c]{3}{*}{40} & 1 & 11 & 7.69 &  & 82.851 & 343.005 & 72.625 &  & 1.1 & $t\ell$ & 1.6 &  & 2.4 & 1.4 &  & 0.20 & 2.75 \\
 &  & 2 & 13 & 6.54 &  & 126.728 & 650.960 & 115.706 &  & 1.3 & $t\ell$ & 1.6 &  & 2.1 & 1.3 &  & 0.23 & 3.32 \\
 &  & 3 & 9 & 11.92 &  & 45.100 & 226.157 & 33.876 &  & 1.2 & $t\ell$ & 1.4 &  & 3.5 & 2.5 &  & 0.14 & 3.29 \\
\end{longtable}
\endgroup

\end{document}